\documentclass[a4paper,12pt]{report}
\usepackage{amsfonts,amssymb,amsmath}
\usepackage{graphicx}
\usepackage{epsfig}
\usepackage[applemac]{inputenc}
\usepackage{latexsym}
\usepackage{stmaryrd}

\newcommand{\HRule}{\rule{\linewidth}{0.5mm}}

\title{Oriented paths in $n$-chromatic digraphs}

\author{Rajai Nasser}

\begin{document}

\begin{titlepage}

\begin{center}


\includegraphics[width=0.15\textwidth]{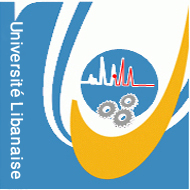}\\[1cm]

\textsc{\LARGE Lebanese University}\\[1.5cm]

\textsc{\LARGE Doctoral Faculty of Sciences and \\[2mm] Technology}\\[1.5cm]

\textsc{\Large Master Thesis}\\[1.5cm]

\HRule \\[0.4cm]
{ \huge \bfseries Oriented paths in $n$-chromatic digraphs}\\[0.4cm]

\HRule \\[2.5cm]

\begin{minipage}{0.4\textwidth}
\begin{flushleft} \large
\emph{Author:}\\
Rajai \textsc{Nasser}
\end{flushleft}
\end{minipage}
\hfill\begin{minipage}{0.25\textwidth}
\emph{Supervisor:} \\
Prof.~Amin \textsc{El-Sahili}
\end{minipage}\\[1.5cm]

\begin{minipage}[center]{0.31\textwidth}
\large
\emph{Jury:}\\
Prof.~Amin \textsc{El-Sahili}\\
Prof.~Hassan \textsc{Abbass}\\
Prof.~Bassam \textsc{Mourad}\\
\end{minipage}

\vfill

{\large July 1, 2010}

\end{center}

\end{titlepage}

\pagenumbering{roman}

\chapter*{Acknowledgement}
\addcontentsline{toc}{chapter}{Acknowledgement}

I am heartily thankful to my supervisor, Prof.~Amin \textsc{El-Sahili}, whose encouragement, guidance and support from the initial to the final level enabled me to develop an understanding of the subject.

\vspace*{4mm}

\noindent I would like to thank Prof.~Hassan \textsc{Abbass} and Prof.~Bassam \textsc{Mourad}, who accepted to be members of the jury.

\vspace*{4mm}

\noindent Lastly, I wish to avail myself of this opportunity, express a sense of gratitude and love to my parents, my brother, my sisters and my friends for their support, help and everything, and to whom I would like to dedicate this work.

\vspace*{4mm}

\noindent Rajai \textsc{Nasser}

\newpage

\chapter*{Abstract}
\addcontentsline{toc}{chapter}{Abstract}

In this thesis, we try to treat the problem of oriented paths in $n$-chromatic digraphs. We first treat the case of antidirected paths in 5-chromatic digraphs, where we explain El-Sahili's theorem and provide an elementary and shorter proof of it. We then treat the case of paths with two blocks in $n$-chromatic digraphs with $n$ greater than 4, where we explain the two different approaches of Addario-Berry et al. and of El-Sahili. We indicate a mistake in Addario-Berry et al.'s proof and provide a correction for it.

\newpage

\tableofcontents

\chapter*{Introduction}
\addcontentsline{toc}{chapter}{Introduction}
\setcounter{page}{1}
\pagenumbering{arabic}
\noindent Gallai-Roy's celebrated theorem ~\cite{pa8, pa9} states that every $n$-chromatic digraph contains a directed path of length $n-1$. More generally, one can ask which connected digraphs are contained in every $n$-chromatic digraph. Such digraphs are called $n$-universal. Since there exist $n$-chromatic graphs with arbitrarily large girth ~\cite{pa10}, $n$-universal digraphs must be oriented trees. Burr ~\cite{pa11} proved that every oriented tree of order $n$ is $(n - 1)^2$-universal (in particular every oriented path is $(n - 1)^2$-universal) and he conjectured that every oriented tree of order $n$ is $(2n - 2)$-universal. This is a generalization of Sumner's conjecture which states that every oriented tree of order $n$ is contained in every tournament of order $2n - 2$. The first linear bound for tournaments was given by H\"{a}ggkvist and Thomason ~\cite{pa12}. The best bound so far, $3n-3$, was obtained by El Sahili ~\cite{pa13}.

\vspace*{4mm}

\noindent Regarding oriented paths in general, there is no better result than the one given by Burr, that is every oriented path is $(n - 1)^2$-universal. However in tournaments, Havet and Thomass\'{e} ~\cite{pa14} proved that except for three particular cases, every tournament of order $n$ contains every oriented path of order $n$.

\vspace*{4mm}

\noindent El-Sahili showed ~\cite{pa1} that except the regular 5-tournament $T_5$, any 5-chromatic oriented digraph in which each vertex has out-degree at least two, contains a copy of the anitidirected path $p_4$ of length 4. To show his result, El-Sahili used a theorem of Gallai ~\cite{pa5}. In chapter II, we give a detailed explanation of El-Sahili's proof, we provide a new elementary shorter proof without using Gallai's theorem, and we conjecture a stronger statement.

\vspace*{4mm}

\noindent El-Sahili conjectured ~\cite{pa15} that every path of order $n \geq 4$ with two blocks is $n$-universal, and Bondy and El-Sahili ~\cite{pa15} proved it if one of the two blocks has length one. The condition $n \geq 4$ is necessary because of odd circuits. El-Sahili and Kouider ~\cite{pa16} introduced the notion of maximal spanning out-forests and used it to show a weak version of El-Sahili's conjecture which states that every path of order $n$ with two blocks is $(n + 1)$-universal.

\vspace*{4mm}

\noindent L. Addario-Berry et al ~\cite{pa2} used strongly connected digraphs and a theorem of Bondy ~\cite{pa17} to show El-Sahili's conjecture. El-Sahili and Kouider ~\cite{pa3} gave a new elementary proof without using strongly connected digraphs or Bondy's theorem. In chapter III we give a detailed explanation of both proofs, we show that there is a small error in Addario-Berry et al' proof and we provide a correction.

\vspace*{4mm}

\noindent All the definitions and basic notations used in this master thesis will be explained in Chapter I.

\chapter{Definitions and basic notations}

\section{Graphs and multi-graphs}

\noindent A \emph{graph} is a pair $G = (V, E)$ of sets such that $E$ is a subset of the power set $P(V)$ of $V$ where every element of $E$ contains exactly two elements of $V$. The elements of $V$ are called the \emph{vertices} of $G$ and the elements of $E$ are called \emph{the edges} of $G$. The set of vertices of $G$ is referred to as $V(G)$, and the set of edges is referred to as $E(G)$. An edge $\{x,y\}$ is noted by $xy$. \emph{The order} $|G|$ of the graph $G$ is the number of vertices in $V(G)$. A graph where we can find an edge between any two distinct vertices is called \emph{complete}. A complete graph of order $n$ is denoted $K_n$;

\vspace*{4mm}

\noindent A \emph{multi-graph} is a triplet $G=(V,E,\varphi)$ where $V$ and $E$ are two sets, and $\varphi$ is a mapping from $E$ into $P(V)$ such that for every $e$ in $E$, $\varphi(e)$ contains one or two vertices of $V$. We say that $V$ is the set of vertices of $G$ and we write $V(G)=V$, similarly we say that $E$ is the set of edges of $G$ and we write $E(G)=E$. \emph{The order of a multi-graph} is also the number of vertices in $V(G)$.

\vspace*{4mm}

\noindent If $e$ is an edge and $\varphi(e)$ contains only one vertex $v$ we say that $e$ is a \emph{loop} on $v$. If $e_1$ and $e_2$ are two different edges on the same vertices i.e. $\varphi(e_1)=\varphi(e_2)$, we say that $e_1$ and $e_2$ are \emph{parallel edges}. A multi-graph $G=(V,E,\varphi)$ without loops or parallel edges can be seen as a graph: we identify it with the graph $G'=(V,\varphi(E))$.

\vspace*{4mm}

\noindent If $G_1$ and $G_2$ are two graphs such that $V(G_1) \subset V(G_2)$ and $E(G_1) \subset E(G_2)$ we say that $G_1$ is a \emph{subgraph} of $G_2$. If in addition $E(G_1)$ contains all the edges $xy$ of $G_2$ such that $x,y \in V(G_1)$, we say that $G_1$ is an \emph{induced subgraph} of $G_2$, and we write $G_1=G_2[V(G_1)]$. If $G_1$ is a subgraph of $G_2$ and $V(G_1)=V(G_2)$ we say that $G_1$ \emph{spans} $G_2$.

\vspace*{4mm}

\noindent A mapping $f: V(G_1) \longrightarrow V(G_2)$ is said to be a \emph{morphism of graphs} if $\forall x,y\in V(G_1)$ we have $f(x)f(y)\in E(G_2)$ whenever $xy\in E(G_1)$. If $f$ is injective, we say that $G_2$ contains a copy of $G_1$ which is $f(G_1):=(f(V(G_1)),\{f(x)f(y)\in E(G_2)/xy\in E(G_1)\})$, or for simplicity we may say that $G_2$ contains $G_1$. If $f$ is bijective, we say that $f$ is an \emph{isomorphism of graphs} and that $G_1$ and $G_2$ are \emph{isomorphic}.

\vspace*{4mm}

\noindent If $G_1=(V_1,E_1,\varphi_1)$ and $G_2=(V_2,E_2,\varphi_2)$ are two multi-graphs, then we say that $G_1$ is a \emph{sub-multi-graph} of $G_2$ if $V_1\subset V_2$, $E_1\subset E_2$ and $\varphi_1$ is the restriction of $\varphi_2$ on $E_1$.

\section{Digraphs and oriented multi-graphs}

\noindent A digraph is a pair $D=(V,E)$ of sets such that $E\subset V\times V$, and such that for every $(x,y)\in E$ we must have $(y,x)\notin E$, in particular if $(x,y)\in E$ then $x\neq y$. We call $V$ the set of vertices of $D$ and we write $V(D)=V$, similarly we call $E$ is the set of arcs (or edges) of $D$ and we write $E(D)=E$. If $e=(x,y)\in E$, we write $x\rightarrow y$; we say that $x$ is the \emph{tail} of $e$ and we write $t(e)=x$ and we say that $y$ is the \emph{head} of $e$ and we write $h(e)=y$. \emph{The order of a digraph} is the number of vertices in $V(D)$.

\vspace*{4mm}

\noindent If $D_1$ and $D_2$ are two digraphs such that $V(D_1) \subset V(D_2)$ and $E(D_1) \subset E(D_2)$ we say that $D_1$ is a \emph{subdigraph} of $D_2$. If in addition $E(D_1)$ contains all the arcs $(x,y)$ of $D_2$ such that $x,y \in V(D_1)$, we say that $D_1$ is an \emph{induced subdigraph} of $D_2$, and we write $D_1=D_2[V(D_1)]$. If $D_1$ is a subdigraph of $D_2$ and $V(D_1)=V(D_2)$ we say that $D_1$ \emph{spans} $D_2$.

\vspace*{4mm}

\noindent A mapping $f: V(D_1) \longrightarrow V(D_2)$ is said to be a \emph{morphism of digraphs} if $\forall x,y\in V(D_1)$ we have $(f(x),f(y))\in E(D_2)$ whenever $(x,y)\in E(D_1)$. If $f$ is injective, we say that $D_2$ contains a copy of $D_1$ which is $f(D_1):=(f(V(D_1)),\{(f(x),f(y))\in E(D_2)/(x,y)\in E(D_1)\})$, or for simplicity we may say that $D_2$ contains $D_1$. If $f$ is bijective, we say that $f$ is an \emph{isomorphism of digraphs} and that $D_1$ and $D_2$ are \emph{isomorphic}.

\vspace*{4mm}

\noindent Let $D=(V,E)$ be a digraph. the \emph{underlying graph} $G(D)$ of $D$ is defined as $G(D):=(V,\psi(E))$, where $\psi: V \times V \longrightarrow P(V)$ is defined as $\psi((x,y))=\{x,y\},\;\forall x,y\in V$. A digraph whose underlying graph is complete is called a \emph{tournament}.

\vspace*{4mm}

\noindent An \emph{Oriented multi-graph} is a triplet $D=(V,E,\varphi)$ where $V$ and $E$ are two sets, and $\varphi$ is a mapping from $E$ into $V \times V$. We say that $V$ is the set of vertices of $D$ and we write $V(D)=V$, similarly we say that $E$ is the set of arcs (or edges) of $D$ and we write $E(D)=E$. If $e\in E(D)$ and $\varphi(e)=(x,y)$ we write $x\rightarrow y$; we say that $x$ is the \emph{tail} of $e$ and we write $t(e)=x$ and we say that $y$ is the \emph{head} of $e$ and we write $h(e)=y$. \emph{The order of an oriented multi-graph} is the number of vertices in $V(D)$.

\vspace*{4mm}

\noindent The \emph{underlying multi-graph} of an oriented multi-graph $D=(V,E,\varphi)$, is the multi-graph $G(D)=(V,E,\psi \circ \varphi)$ where $\psi: V \times V \longrightarrow P(V)$ is defined as $\psi((x,y))=\{x,y\},\;\forall x,y\in V$. If $D$ is an oriented multi-graph whose underlying multi-graph is a graph (contains no loops and no parallel edges), $D$ can be seen as a digraph: we identify $D$ with $D'=(V,\varphi(E))$.

\vspace*{4mm}

\noindent If $D_1=(V_1,E_1,\varphi_1)$ and $D_2=(V_2,E_2,\varphi_2)$ are two oriented multi-graphs, then $D_1$ is a \emph{sub-oriented-multi-graph} of $D_2$ if $V_1\subset V_2$, $E_1\subset E_2$ and $\varphi_1$ is the restriction of $\varphi_2$ on $E_1$.

\vspace*{4mm}

\noindent For simplicity, we will not be strict when dealing with graphs (resp. multi-graphs, oriented multi-graphs or digraphs), in the sense that if $G$ is a graph (resp. multi-graph, oriented multi-graph or digraph) we may not differ strictly between $G$ and $V(G)$ or between $G$ and $E(G)$: If $v$ is a vertex of $G$ and $e$ is an edge of $G$, we may write $v\in G$ and $e\in G$ rather than $v\in V(G)$ and $e\in E(G)$. Also if $H$ is a subgraph (resp. sub-multi-graphs, sub-oriented-multi-graph or subdigraph) of $G$, and $e$ is an edge (or arc) of $G$ we denote by $H\cup e$ or $H+e$ the subgraph (resp. sub-multi-graphs, sub-oriented-multi-graph or subdigraph) of $G$ obtained from $H$ by adding the edge $e$, and if $e\in H$ we denote by $H-e$ the subgraph (resp. sub-multi-graphs, sub-oriented-multi-graph or subdigraph) of $G$ obtained from $H$ by deleting the edge $e$.

\section{Degree and neighborhood of a vertex}

\noindent Let $G$ be a graph, if $e=xy$ is an edge of $G$ and we say that the vertices $x$ and $y$ are \emph{adjacent} and we say that $e$ is \emph{incident} to $x$ and $y$. \emph{The neighborhood} $N(v)$ of a vertex $v$ is defined as the set of vertices adjacent to it, and its \emph{degree} $d(v)$ is the number of vertices in $N(v)$ which is equal to the number of edges incident to $v$.

\vspace*{4mm}

\noindent Let $G$ be a graph. The \emph{maximum degree} of $G$ is defined as $\Delta(G):=max\{d(v)/v\in V(G)\}$ and the \emph{minimum degree} of $G$ is defined as $\delta(G):=min\{d(v)/v\in V(G)\}$.

\vspace*{4mm}

\noindent Let $D=(V,E)$ be a digraph. \emph{The neighborhood} $N(v)$ of a vertex $v$ is its neighborhood in the underlying graph. \emph{The degree} $d(v)$ of a vertex $v$ is its degree in the underlying graph. the \emph{out-neighborhood} $N^+(v)$ of a vertex $v$ is defined as $N^+(v):=\{w\in V(D)/v\rightarrow w\}$. Similarly the \emph{in-neighborhood} $N^-(v)$ is defined as $N^-(v):=\{w\in V(D)/w\rightarrow v\}$. The \emph{out-degree} $d^+(v)$ of a vertex $v$ is the number of arcs whose tail is $v$, and the \emph{in-degree} $d^-(v)$ of $v$ is the number of arcs whose head is $v$. We define also $\Delta^+(D):=max\{d^+(v)/v\in V(D)\}$, $\delta^+(D):=min\{d^+(v)/v\in V(D)\}$, $\Delta^-(D):=max\{d^-(v)/v\in V(D)\}$ and $\delta^-(D):=min\{d^-(v)/v\in V(D)\}$.

\section{Paths and cycles}

\noindent Let $G$ be a graph (or multi-graph), a \emph{path} $P$ from $x$ to $y$ in $G$ is a finite sequence $P=x_1x_2...x_n$ (we can have $n=1$) of distinct vertices such that $x_1=x$, $x_n=y$ and $x_ix_{i+1} \in E(G)$ for $1\leq i\leq n-1$, $x_1$ and $x_n$ are the \emph{end vertices} of $P$. A \emph{sub-path} of a path $P$ is a path which is a subset of $P$. A \emph{cycle} $C$ in $G$ is a finite sequence $C=x_1x_2...x_n$ (we can have $n=1$) of distinct vertices such that $x_nx_1\in E(G)$ and $x_ix_{i+1} \in E(G)$ for $1\leq i\leq n-1$. If $C=x_1...x_n$ is a cycle, then if $e=x_ix_j \in E(G)$ such that $i-j\neq 1$ mod $n$ and $j-i\neq 1$ mod $n$, then $e$ is called a \emph{chord} of $C$, if such chord does not exist we say that $C$ is \emph{chordless}. A graph (or multi-graph) is said to be \emph{acyclic} if it does not contain any cycle; note that an acyclic multi-graph is necessarily a graph. \emph{The length} $l(P)$ (resp. $l(C)$) of a path $P$ (resp. a cycle $C$) is the number of edges in it. A \emph{hamiltonian} path $P$ (resp. cycle $C$) is a path (resp. cycle) which spans $G$, i.e. $V(P)=V(G)$ (resp. $V(C)=V(G)$).

\vspace*{4mm}

\noindent \emph{The distance} $d(x,y)$ between two vertices $x$ and $y$, is the minimal length of a path from $x$ to $y$ if such path exists. If there is no path between $x$ and $y$, we set $d(x,y):=\infty$. The map $d: V(G)\times V(G) \longrightarrow \mathbb{R}\cup \{\infty\}$ verifies the axioms of generalized metric, and so $(V(G),d)$ is a generalized metric space. \emph{The diameter} $d(G)$ of a graph (or multi-graph) $G$, is the maximal distance between two vertices of $G$.

\vspace*{4mm}

\noindent \emph{The girth} $g(G)$ of a graph (or multi-graph) $G$ is the minimal length of a cycle in it if such one exists, and if $G$ is acyclic we set $g(G):=\infty$. If $g(G)=1$ then $G$ contains a loop and if $g(G)=2$ then $G$ contains parallel edges, so if $g(G)\geq3$, $G$ contains no loops and no parallel edges and so $G$ is necessarily a graph. Note that all the above definitions are also defined for oriented multi-graphs by applying them on the underlying multi-graphs.

\vspace*{4mm}

\noindent Let $D$ be a digraph. A \emph{directed path} $P$ in $D$ is a finite sequence of different vertices $P=x_1x_2...x_n$ (we can have n=1) such that $x_i\rightarrow x_{i+1}$ for $1\leq i\leq n-1$. A \emph{block of a path} $P$ in $D$ is a maximal directed sub-path of $P$. A path having $l$ blocks of consecutive lengths $k_1,k_2,...,k_l$ is denoted by $P^+(k_1,k_2,...,k_l)$ (or $P(k_1,k_2,...,k_l)$) if $x_1\rightarrow x_2$ and $P^-(k_1,k_2,...,k_l)$ if $x_1\leftarrow x_2$. An \emph{antidirected} path is a path whose blocks are all of length 1. A \emph{circuit} $C$ in $D$ is a finite sequence of different vertices $C=x_1x_2...x_n$ (we can have n=1) such that $x_i\rightarrow x_{i+1}$ for $1\leq i\leq n-1$ and $x_n\rightarrow x_1$.

\section{Connectivity}

\noindent Let $G$ be a graph (or multi-graph), $G$ is \emph{connected} if $d(G)<\infty$, i.e. there exist a path between any two vertices. $G$ is \emph{disconnected} if it is not connected. $G$ is \emph{$k$-connected} if it remains connected after the removal of any $k'<k$ vertices. The \emph{connectivity} $\kappa(G)$ is the maximal integer $k$ such that $G$ is $k$-connected ($\kappa(G)=0$ if and only if $G$ is disconnected). $G$ is \emph{$k$-edge-connected} if it remains connected after the removal of any $k'<k$ edges. The \emph{edge-connectivity} $\lambda(G)$ is the maximal integer $k$ such that $G$ is $k$-edge-connected ($\lambda(G)=0$ if and only if $G$ is disconnected).

\vspace*{4mm}

\noindent Let $G$ be a graph (or multi-graph), a maximal connected subgraph of $G$ is called a \emph{connected component} of $G$. Suppose that $G$ is connected, a vertex $v$ whose removal disconnect $G$ is a \emph{cut-vertex} of $G$ and an edge $e$ whose removal disconnect $G$ is a \emph{bridge} of $G$. A maximal connected subgraph of $G$ without cut-vertices is called a \emph{block} of $G$.

\vspace*{4mm}

\noindent Let $D$ be a digraph (resp. oriented multi-graph), all the above notations are defined for $D$ by applying them on its underlying graph (resp. underlying multi-graph). $D$ is called \emph{strongly connected} if by choosing any two vertices $x$ and $y$ in $D$ we can find a directed path from $x$ to $y$ and a directed path from $y$ to $x$.

\section{Trees and forests}

\noindent An acyclic graph is called a \emph{forest}. A connected acyclic graph is called a \emph{tree}, so a forest is the union of trees (each graph is the union of its connected components). The vertices of degree 1 in a tree are called the \emph{leaves} of the tree. A tree containing only one vertex is called a \emph{trivial tree}, then a non-trivial tree contain at least two leaves (consider for example the ends of a longest path).

\vspace*{4mm}

\noindent The following assertions are equivalent for a graph T (the proof is straightforward for the first four, use a simple induction for the last two):
\begin{enumerate}
  \item $T$ is a tree.
  \item Any two vertices of $T$ are linked by a unique path.
  \item $T$ is minimally connected, i.e. $T$ is connected and $T-e$ is disconnected for all edges $e\in E(T)$.
  \item $T$ is maximally acyclic, i.e $T$ is acyclic and $T+xy$ contains a cycle for any two non-adjacent vertices $x$ and $y$ of $G$.
  \item $T$ is connected and $|E(T)|=|V(T)|-1$.
  \item $T$ is acyclic and $|E(T)|=|V(T)|-1$.
\end{enumerate}

\vspace*{4mm}

\noindent An \emph{oriented tree} is a digraph whose underlying graph is a tree, similarly an \emph{oriented forest} is a digraph whose underlying graph is a forest. An \emph{out-leaf} of a tree is a leaf whose out-degree is zero, similarly an \emph{in-leaf} of a tree is a leaf whose in-degree is zero. An \emph{out-branching} (resp. \emph{in-branching}) is an oriented tree in which a unique vertex which we call \emph{the root} has its in-degree (resp. out-degree) 0, and the other vertices has in-degree (resp. out-degree) 1. An \emph{out-forest} (resp. \emph{in-forest}) is an oriented forest whose connected components are out-branchings (resp. in-branchings). Let $F$ be an out-forest, the level $l_F(v)$ of a vertex $v\in F$ is the order of a longest directed path ending at $v$.

\section{Coloring}

\noindent A \emph{$k$-coloring} of a graph (or multi-graph) $G$ is a mapping $c:G\longrightarrow \{1,2,...,k\}$ (we can use any set of $k$ elements instead of $\{1,2,...,k\}$). If $v$ is a vertex of $G$ we say that $c(v)$ is the color of $v$, and if $v$ is adjacent to a vertex of color $i$, we say that $v$ is adjacent to the color $i$. A \emph{good $k$-coloring} of a graph $G$ is a coloring $c$ such that any adjacent vertices does not have the same color.

\vspace*{4mm}

\noindent If $G$ admits a good $k$-coloring, we say that $G$ is \emph{$k$-colorable}. A subset $L$ of $V(G)$ is said to be \emph{stable} if there is no adjacent vertices in it, i.e. the set of edges in the subgraph $G[L]$ of $G$ induced by $L$ is empty. $G$ is said to be \emph{independent} if $V(G)$ is stable. Note that $G$ is $k$-colorable if and only if we can partition $V(G)$ into $k$ stable subsets.

\vspace*{4mm}

\noindent The \emph{chromatic number} $\chi(G)$ of $G$ is the least integer $k$ such that $G$ is $k$-colorable. If $\chi(G)=k$ and $\chi(G-v)<k\;\forall v\in V(G)$ we say that $G$ is \emph{$k$-critical}. All the above notations can be defined for digraphs (resp. oriented multi-graphs) by applying them on their underlying graphs (resp. underlying multi-graphs).

\section{Contraction and minors}

\noindent Let $G$ be a graph, and let $H$ be a subset of $V(G)$ (or a subgraph of G), then the graph obtained from $G$ by \emph{contracting} $H$ is $G/H$ defined by $V(G/H):=(V(G)\setminus H)\cup\{v_H\}$ where $v_H$ is a new vertex and $E(G/H):=\{xy / xy\in E(G),\;x,y\in V(G)\setminus H\}\cup\{vv_H/v\in V(G)\setminus H,\;\exists v'\in H,\;vv'\in E(G)\}$.

\vspace*{4mm}

\noindent Let $D=(V,E)$ be a digraph, and let $H$ be a subdigraph of $D$, We say that $D$ is \emph{contractable} by $H$ if for all vertices $v$ in  $V(D)\setminus V(H)$, we cannot find two arcs $v\rightarrow x$ and $y\rightarrow v$ such that $x\in H$ and $y\in H$, i.e. all arcs in $D$ between $v$ and $H$ are in the same direction. In this case, the digraph obtained from $D$ by \emph{contracting} $H$ is $D/H$ defined by $V(D/H):=(V(D)\setminus H)\cup\{v_H\}$ where $v_H$ is a new vertex and $E(D/H):=\{(x,y) / (x,y)\in E(D),\;x,y\in V(D)\setminus H\}\cup\{(v_H,v)/v\in V(G)\setminus H,\;\exists v'\in H,\;(v',v)\in E(G)\}\cup\{(v,v_H)/v\in V(G)\setminus H,\;\exists v'\in H,\;(v,v')\in E(G)\}$.

\vspace*{4mm}

\noindent Let $D=(V,E,\varphi)$ be an oriented multi-graph, and let $H$ be a sub-oriented-multi-graph of $D$, then the oriented multi-graph obtained from $D$ by \emph{contracting} $H$ is $D/H=(V',E',\varphi')$ where $V'=(V\setminus V(H))\cup\{v_H\}$, $E'=E\setminus E(H)$ and $\varphi':E' \longrightarrow P(V')$ is defined by $\varphi'(e)=(f_H(t_D(e)),f_H(h_D(e)))\;\forall e\in E'$ where $f_H(x)=x$ if $x\notin H$ and $f_H(x)=v_H$ if $x\in H$.

\vspace*{4mm}

\noindent Note that the notation $D/H$ have different meaning when interpreting $D$ as a digraph or oriented multi-graph. The notation takes its meaning relatively to the context.

\vspace*{4mm}

\noindent If $G$ is a graph (resp. digraph or oriented multi-graph) and if $G'$ is a graph (resp. digraph or oriented multi-graph), we say that $G'$ is a \emph{minor} of $G$ if there exist a finite sequence $G_1,G_2,...,G_n$ of graphs (resp. digraph or oriented multi-graph) such that $G_1=G$, $G_n=G'$ and $\forall i\in \{1,2,...,n-1\}$ $G_{i+1}$ is a subgraph (resp. subdigraph or sub-oriented-multi-graph) of $G_i$ or obtained from $G_i$ by contracting some subgraph (resp. subdigraph or sub-oriented-multi-graph) of it.

\chapter{Antidirected paths in digraphs}

\section{Introduction}

\noindent \emph{The antidirected path} $p_4$ is a digraph defined, up to isomorphism as follows:

\vspace*{1mm}

\noindent $$V(p_4)=\{x,y,z,v,w\}, E(p_4)=\{(y,x),(y,z),(v,z),(v,w)\}$$

\vspace*{2mm}

\noindent Let $T_5$ be the 5-tournament satisfying $d^+(u)=d^-(u)=2 \;\forall u \in T_5$. Grunbaum ~\cite{pa5} proved that $T_5$ is the only 5-tournament which doesn't contain a copy of $p_4$. El-Sahili ~\cite{pa1} showed that except $T_5$, any 5-chromatic oriented digraph in which each vertex has out-degree at least two, contains a copy of $p_4$. He showed by an example that the condition that each vertex has out-degree at least two is necessary.

\vspace*{4mm}

\noindent To show his result, El-Sahili used a theorem of Gallai ~\cite{pa6}, which states that if $G$ is k-critical, then each block of the subgraph of $G$ induced by the vertices of degree $k-1$, is either complete or chordless odd cycle.

\vspace*{4mm}

\noindent In this chapter we will give a detailed explanation of the argument used by El-Sahili to show his theorem. We will then provide a new elementary shorter proof which does not require the use of Gallai's theorem. We conclude this chapter by stating a new conjecture generalizing this theorem.

\section{First Step of the proof}

\noindent {\bf Theorem 2.1} ~\cite{pa1}: \emph{Let $D$ be a 5-chromatic connected digraph distinct from $T_5$ in which each vertex has out-degree at least two. Then $D$ contains a copy of $p_4$}.

\vspace*{6mm}

\noindent To prove this theorem, we need several lemmas:

\vspace*{6mm}

\noindent {\bf Lemma 2.2} ~\cite{pa6}: \emph{Except for $T_5$, any 5-tournament contains a copy of $p_4$}.

\vspace*{2mm}

\noindent {\bf Proof}: Suppose to the contrary that there exists a 5-tournament $T$ other then $T_5$ which does not contain any $p_4$, then there exist at least one vertex $v_1$ in $T$ such that $d^+(v)\geq3$, let $\{v_2,v_3,v_4\}\subset N^+(v)$. we can assume without loss of generality that we have $v_2 \rightarrow v_3 \rightarrow v_4 \rightarrow v_2$ or $v_2 \rightarrow v_3 \rightarrow v_4 \leftarrow v_2$.

\vspace*{4mm}

\noindent In the first case, if $\exists i \in \{2,3,4\}$ such that $v_i\rightarrow v_5$, we may assume without loss of generality that $i=2$, then $v_5v_2v_3v_1v_4$ is a $p_4$, so we conclude that $v_5\rightarrow v_i, \forall i\in\{2,3,4\}$, but in this case $v_3v_5v_2v_1v_4$ would be a $p_4$.

\vspace*{4mm}

\noindent In the latter case, we have $v_5\rightarrow v_2$ because otherwise $v_5v_2v_3v_1v_4$ is a $p_4$. If $\exists i \in \{3,4\}$ such that $v_5 \rightarrow v_i$, we may assume without loss of generality that $i=3$ and so $v_3v_5v_2v_1v_4$ is be a $p_4$, so we conclude that $\forall i\in\{3,4\}, v_i\rightarrow v_5$, but in this case $v_5v_3v_4v_1v_2$ is a $p_4$. \hfill $\boxempty$

\vspace*{6mm}

\noindent {\bf Corollary 2.3}: \emph{If $D$ is a digraph verifying the conditions of Theorem 2.1 and if $D$ contains a $K_5$, then $D$ contains a copy of $p_4$}.

\vspace*{2mm}

\noindent {\bf Proof}: Let $p_3$ be the subpath of $p_4$ formed by the first three edges, and let $p_2$ be the subpath of $p_4$ formed by the first two edges. Since $G(D)$ contains $K_5$, then $D$ contains a 5-tournament. If this 5-tournament is not $T_5$ then by \emph{Lemma 2.2} we conclude that $D$ contains a copy of $p_4$.

\vspace*{4mm}

\noindent Then we may assume that $D$ contains $T_5$, and since $D$ is not exactly $T_5$, then we will have an edge $xy$ in $G(D)$ such that $x$ is outside $T_5$ and $y$ belongs to $T_5$. If $y \rightarrow x$ then this edge along with a path $p_3$ in $T_5$ starting at $y$ (we can always find a copy of $p_3$ in $T_5$ starting at any point of it), form a path $p_4$. Otherwise, since $d^+(x)\geq2$, then there exist a vertex $z$ distinct from $y$ such that $x \rightarrow z$. If $z \notin T_5$, then the path $zxy$ along with a copy of $p_2$ in $T_5$ starting at $y$, form a copy of $p_4$. If $z \in T_5$, then the path $zxy$ along with a copy of $p_2$ in $T_5$ starting at $y$ and not intersecting $z$ (For any two vertices of $T_5$, we can always find a copy of $p_2$ starting at one vertex and not intersecting the other), form a copy of $p_4$. \hfill $\boxempty$

\vspace*{6mm}

\noindent {\bf Theorem 2.4} ~\cite{pa7}: \emph{If $G$ is a connected graph which is not complete nor an odd cycle, then $\chi(G)\leq \Delta(G)$}.

\vspace*{6mm}

\noindent {\bf Corollary 2.5} ~\cite{pa4}: \emph{If $D$ is a digraph which does not contain any tournament of order $2n+1$ ($n\geq2$), and in which any vertex has in-degree at most $n$, then $\chi(D)\leq2n$}.

\vspace*{2mm}

\noindent {\bf Proof}: Suppose to the contrary that $\chi(D)\geq2n+1$ and let $D'$ be a $2n+1$-critical subdigraph of $D$. If there exists a vertex $v$ in $D'$ such that $d_{D'}^+(v)<n$ then $d(v)<2n$, and since $D'$ is $2n+1$-critical then $\chi(D'-v)=2n$ and we can easily check that $\chi(D')=2n$ (extend a good $2n$-coloring of $D'-v$ by giving $v$ a color not adjacent to it; we can find such color since $v$ is adjacent to at most $2n-1$ vertices) which contradicts the fact that $\chi(D')=2n+1$.

\vspace*{4mm}

\noindent We conclude that for every vertex in $D'$ we have $d_{D'}^+(v)\geq n$ and $d_{D'}^-(v)\leq n$, and since $\sum_{v \in D'} d_{D'}^+(v) = \sum_{v \in D'} d_{D'}^-(v) = |E(D')|$ we conclude that for every vertex $v$ in $D'$ we have $d_{D'}^+(v) = d_{D'}^-(v) = n$ and so $d_{D'}(v)=2n$ which implies that $\Delta(D')=2n$. Obviously $G(D')$ is not an odd cycle since $\Delta(D')=2n\geq2\times2=4$, and it is not complete since otherwise $D$ would contain a $2n+1$-tournament, so by Brooks theorem (Theorem 2.4) we conclude that $\chi(D')\leq\Delta(D')=2n$ which contradicts the fact that $\chi(D')=2n+1$. \hfill $\boxempty$

\vspace*{6mm}

\noindent {\bf Lemma 2.6}: \emph{If $D$ is a connected digraph in which each vertex has in-degree at most one. Then $D$ contains at most one cycle which is a circuit}.

\vspace*{2mm}

\noindent {\bf Proof}: Let $C$ be a cycle which is a subdigraph of $D$, if $C$ is not a circuit then $\exists v \in C$ such that $d_C^-(v)\neq1$ or $d_C^+(v)\neq1$. Since $d_C(v)=2$ and $d_C^-(v)\leq1$ then we have $d_C^-(v)=0$ and $d_C^+(v)=2$, but we have $\displaystyle\sum_{v \in C}d_C^+(v)=\sum_{v \in C}d_C^-(v)=n$ then  $\exists w \in C$ such that $d_C^-(w)=2$ which is a contradiction. We conclude that every cycle in $D$ is necessarily a circuit.

\vspace*{4mm}

\noindent Suppose that there exist two different circuits $C_1$ and $C_2$ subdigraphs of $D$. Suppose that $C_1 \cap C_2 \neq \phi$, since $C_1 \neq C_2$ then we can say without loss of generality that $C_2 \nsubseteq C_1$ and thus $\exists v \in C_2 \setminus C_1$ such that $\exists w \in C_1$ with $v \rightarrow w$, but $w$ has another in-neighbor in $C_1$ which is a contradiction. We conclude that $C_1 \cap C_2 = \phi$, but $D$ is connected then there exists a path between a vertex of $C_1$ and a vertex of $C_2$, let $P=x_1x_2...x_n$ be a minimal such path ($x_1 \in C_1$ and $x_n \in C_2$). $P$ is minimal, so $x_1$ is the only vertex of $P$ in $C_1$ and $x_2$ is the only vertex of $P$ in $C_2$. Since $x_1$ has an in-neighbor in $C_1$ and $d^-(x_1)\leq1$ we have $x_1 \rightarrow x_2$, let $i$ be the maximum integer such that $x_i \rightarrow x_{i+1}$. If $i<n-1$ then $x_i \rightarrow x_{i+1}$ and $x_{i+2} \rightarrow x_{i+1}$ which contradicts the fact that $d^-(x_{i+1})\leq1$, so $i=n-1$ and $x_{n-1} \rightarrow x_n$ but $x_n$ has another in-neighbor in $C_2$ which contradicts the fact that $d^-(x_n)\leq1$. \hfill $\boxempty$

\vspace*{6mm}

\noindent Note that the above corollary and lemma holds also when we substitute ``in-degree" by ``out-degree".

\vspace*{6mm}

\noindent In the sequel, $D$ denote an oriented digraph verifying the conditions of \emph{theorem 2.1}. We suppose to the contrary that $D$ does not contain any copy of $p_4$. by the above corollary we may assume that $D$ does not contain any 5-tournament. Let $D'$ be a 5-critical subdigraph of $D$ and let $D^o$ be the subdigraph of $D'$ induced by the vertices of out-degree at least three in $D'$, i.e. $D^o=\{x \in D'/d_{D'}^+ \geq 3\}$.

\vspace*{6mm}

\noindent {\bf Lemma 2.7}: \emph{$D'$ contains at least one vertex whose out-degree in $D'$ is at least three, i.e. $D^o$ is not empty}.

\vspace*{2mm}

\noindent {\bf Proof}: Otherwise we would have $d_{D'}^+(v)\leq2$ for every vertex $v$ in $D'$. $D$, and hence $D'$, contains no 5-tournament, so by \emph{corollary 2.5} we conclude that $\chi(D')\leq4$, which is a contradiction. \hfill $\boxempty$

\vspace*{6mm}

\noindent {\bf Lemma 2.8}: \emph{Every vertex in $D$ has at most one in-neighbor in $D^o$}.

\vspace*{2mm}

\noindent {\bf Proof}: Suppose to the contrary that there exists a vertex $v$ having two in-neighbors $x,y \in D^o$ and let $\{v,x_1,x_2\} \subset N^+(x)$. If $y \in \{x_1,x_2\}$, we may suppose without loss of generality that $y=x_1$, then $\exists y_1\in N^+(y)\setminus\{v,x,x_1,x_2\}$ since $d_{D'}^+(y)\geq3$, thus $y_1yvxx_2$ is a $p_4$, a contradiction. So $y \notin \{x_1,x_2\}$ and more generally we can say that $x$ and $y$ are not adjacent. $d^+(y)\geq3$ so $\exists y_1 \in N^+(y)\setminus\{x_1,x,v\}$ thus $x_1xvyy_1$ is a $p_4$, a contradiction. \hfill $\boxempty$

\vspace*{6mm}

\noindent {\bf Corollary 2.9}: \emph{$\forall v\in D^o, d_{D^o}^-(v)\leq1$}.

\vspace*{2mm}

\noindent {\bf Proof}: Clear. \hfill $\boxempty$

\vspace*{6mm}

\noindent {\bf Lemma 2.10}: \emph{Let $v$ be a vertex of $D$ such that $d^+(v)\geq3$ and $\{x,y,z\} \subset N^+(v)$. If $x \rightarrow y$ then $x\rightarrow z$, $yz \notin E(G(D))$ and $N^-(y)=N^-(z)=\{v,x\}$}.

\vspace*{2mm}

\noindent {\bf Proof}: If $x \nrightarrow z$ then $\exists w \in N^+(x)\setminus \{v,y,z\}$ since $d^+(v)\geq2$, so $wxyvz$ is a $p_4$ which is a contradiction. So we must have $x\rightarrow z$.

\noindent If we suppose that $yz \in E(G(D))$, we may assume without loss of generality that $y \rightarrow z$. We have $d^+(y)\geq2$ so $\exists w \in N^+(y)\setminus\{v,x,y,z\}$ and then $wyzvx$ is a $p_4$, a contradiction. So we have $yz \notin E(G(D))$.

\noindent Suppose that $N^-(y) \neq \{v,x\}$, so $\exists w \in N^-(y)\setminus\{v,x,y,z\}$, and since $d^+(w)\geq2$ then $\exists w' \neq y$ such that $w\rightarrow w'$. If $w'=v$ then $vwyxz$ is a $p_4$, a contradiction. So $w'\neq v$, let $u\in\{x,z\}\setminus \{w'\}$, then $w'wyvu$ is a $p_4$, a contradiction. So $N^-(y) = \{v,x\}$ (We prove similarly that $N^-(z) = \{v,x\}$). \hfill $\boxempty$

\vspace*{6mm}

\noindent {\bf Lemma 2.11}: \emph{If $v$ and $v'$ are two vertices such that there exist two adjacent vertices $x$ and $y$ in $N^+(v)\cap N^-(v')$, then $N^+(v)=\{x,y\}$}.

\vspace*{2mm}

\noindent {\bf Proof}: We may assume without loss of generality that $x \rightarrow y$. If $N^+(v)\neq\{x,y\}$ then $\exists w \in N^+(v)\setminus\{x,y\}$, by \emph{lemma 2.10} $w$ cannot be $v'$, and so $wvyxv'$ is a $p_4$ which is a contradiction. \hfill $\boxempty$

\vspace*{6mm}

\noindent {\bf Lemma 2.12}: \emph{$D^o$ is an independent set of $D$}.

\vspace*{2mm}

\noindent {\bf Proof}: Suppose to the contrary that $D^o$ is not an independent set, so there exist a connected component $L$ of $D^o$ which contains at least two vertices. If $L$ is a circuit, then every vertex of $L$ has one in-neighbor in $L$ and has at least two out-neighbors outside $D^o$ since its out-degree in $D'$ is at least 3. If $L$ is not a cycle, let $v$ be the last vertex in a maximal directed path in $L$, we can easily verify that $d_{D^o}^-(v)=1$ and $d_{D^o}^+(v)=0$ so $v$ has at least three out-neighbors outside $D^o$ since $d_{D'}^+(v)\geq3$.

\vspace*{4mm}

\noindent So in all cases, we can always find a vertex $v$ in $L$ having at least two out-neighbors outside $D^o$ and such that $d_L^-(v)=1$, let $v'$ be the in-neighbor of $v$ in $L$ and let $v_1$, $v_2$ and $v_3$ three out-neighbors of $v$ in $D'$ such that $v_1$ and $v_2$ are outside $D^o$ (i.e. $d_{D'}^+(v_1)\leq2$ and $d_{D'}^+(v_2)\leq2$). $D'$ is 5-critical, so any vertex in $D'$ has at least 4 neighbors, we conclude that $d_{D'}^-(v_1)\geq2$ and $d_{D'}^-(v_2)\geq2$. If $v_1$ and $v_2$ are not adjacent, $v_1$ has one in-neighbor in $D'\setminus\{v,v_1,v_2\}$, otherwise we may assume without loss of generality that $v_1\rightarrow v_2$, but since $d_{D'}^-(v_1)\geq2$ we conclude again that $v_1$ has one in-neighbor in $D'\setminus\{v,v_1,v_2\}$. So in all cases we can say without loss of generality that there exist a vertex $u$ in $D'\setminus\{v,v_1,v_2\}$ such that $u \rightarrow v_1$.

\vspace*{4mm}

\noindent Suppose that $u=v_3$, by \emph{Lemma 2.10} we have $v_3 \rightarrow v_2$ and by \emph{Lemma 2.8} we have $u \notin D^o$ so $d_{D'}^-(u)\geq2$ which implies that $\exists w \in D'\setminus\{v,v_1,v_2,v_3\}$ such that $w \rightarrow u$. Since $d^+(w)\geq2$, there exists a vertex $w'\neq u$ such that $w \rightarrow w'$, suppose that $w'\neq v$ then let $w'' \in \{v_1,v_2\}\setminus\{w'\}$ so $w'wuvw''$ is a $p_4$, a contradiction. So we have $w\rightarrow v$. $d^+_{D'}(v')\geq3$ then let $\{v'_1,v'_2,v\}\subset N^+_{D'}(v')$, so by \emph{lemma 2.8} $w\neq v'$ and $\{v'_1,v'_2\}\cap\{v_1,v_2,v_3\}=\phi$. Let $w''\in\{v'_1,v'_2\}\setminus\{w\}$, then $w''v'vwu$ is a $p_4$, a contradiction.

\vspace*{4mm}

\noindent We conclude that $u \notin \{v,v_1,v_2,v_3\}$ and by \emph{lemma 2.8} $u$ cannot be $v'$. $d^+(u)\geq2$ so there exist a vertex $u'$ different from $v_1$ such that $u\rightarrow u'$. If $u'\neq v$, let $w \in \{v_2,v_3\}\setminus\{u'\}$, so $u'uv_1vw$ is a $p_4$, a contradiction. So $u'=v$, and by \emph{lemma 2.11} we cannot have $v' \rightarrow u$, and since $d^+(v')\geq3$, $\exists w \in N^+(v')\setminus\{u,v_1,v\}$ so $wv'vuv_1$ is a $p_4$, a contradiction. \hfill $\boxempty$

\vspace*{6mm}

\noindent {\bf Lemma 2.13}: \emph{Let $v \in D^o$, then $v$ has exactly three out-neighbors $v_1,v_2$ and $v_3$ in $D'$ such that $v_1\rightarrow v_2$ and $v_1\rightarrow v_3$}.

\vspace*{2mm}

\noindent {\bf Proof}: Suppose that any two out-neighbors of $v$ in $D'$ are not adjacent, and let $v_1$, $v_2$ and $v_3$ be three out-neighbors of $v$ in $D'$. $\forall i \in \{1,2,3\}, \exists u_i \in D'\setminus\{v,v_1,v_2,v_3\}$ such that $u_i \rightarrow v_i$. By \emph{Lemma 2.8} we cannot have $u_1=u_2=u_3$ because otherwise we would have $u\in D^o$, so we may assume without loss of generality that $u_1 \neq u_2$. Suppose $u_1 \nrightarrow v$, we have $d^+(u_1)\geq2$ so $\exists w \in D\setminus\{v_1,u_1,v\}$, let $w' \in \{v_2,v_3\}\setminus\{w\}$ so $wu_1v_1vw'$ is a $p_4$ which is a contradiction. So we conclude that $u_1 \rightarrow v$ and similarly $u_2 \rightarrow v$ but we will have a copy of $p_4$ which is $v_1u_1vu_2v_2$, a contradiction. So we conclude that at least two in-neighbors of $v$, say $v_1$ and $v_2$, are adjacent, we may suppose that $v_1\rightarrow v_2$, so by \emph{lemma 2.10} we have also $v_1\rightarrow v_3$.

\vspace*{4mm}

\noindent Suppose that $v$ has four out-neighbors $v_1$,$v_2$,$v_3$ and $v_4$ in $D'$. By \emph{lemma 2.10} we have $v_1 \rightarrow v_3$ and $v_1 \rightarrow v_4$, thus $v_1 \in D^o$ and $\{v,v_1\} \subset N^-(v_2)$ which gives a contradiction with \emph{lemma 2.8}. \hfill $\boxempty$

\section{Second Step (El-Sahili's proof)}

\noindent In the sequel, we will need the use of the following theorem proved by Gallai:

\vspace*{6mm}

\noindent {\bf Theorem 2.14} ~\cite{pa5}: \emph{Let $G$ be a $k$-critical graph, then each block of the subgraph of $G$ induced by the vertices of degree $k-1$ is either complete or a chordless odd cycle}.

\vspace*{6mm}

\noindent Let $D_m$ be the subdigraph of $D'$ induced by the vertices of degree 4.

\vspace*{6mm}

\noindent {\bf Lemma 2.15}: \emph{Any vertex $v$ of $D'$ has at least two in-neighbors in $D'$}.

\vspace*{2mm}

\noindent {\bf Proof}: If $v \in D'\setminus D^o$, then $d_{D'}(v)\geq4$ because $D'$ is 5-critical, and since $d_{D'}^+(v)\leq2$ then $d_{D'}^-(v)\geq2$. So we may assume that $v \in D^o$, let $N_{D'}^+(v)=\{v_1,v_2,v_3\}$ where $v_1 \rightarrow v_2$, $v_1 \rightarrow v_3$, $v_2v_3 \notin E(G(D))$ and $N^-(v_2)=N^-(v_3)=\{v,v_1\}$ (By \emph{lemmas 2.10 and 2.13}). $\forall w \in N_{D'}^-(v_1)\setminus\{v\}$, we have $w \rightarrow v$ because otherwise we can find $u \in N^+(w)\setminus\{v,v_1\}$ and we can find $w' \in \{v_2,v_3\}\setminus\{u\}$, so $w'vv_1wu$ is $p_4$, which is a contradiction. If we suppose that $d_{D'}^-(v)=1$ we will have $d_{D'}^-(v_1)=2$ and so $d_{D'}(v)=d_{D'}(v_1)=d_{D'}(v_2)=d_{D'}(v_3)=4$. Thus $v_1$, $v_2$, $v_3$ and $v_3$ are in the same block of $D_m$, this block cannot be an odd cycle, so by \emph{theorem 2.14} it's complete which contradicts the fact that $v_2v_3 \notin E(G(D))$. \hfill $\boxempty$

\vspace*{6mm}

\noindent We now associate to each vertex $v$ in $D^o$ the set $S(v)=\{t(v),t'(v),v_0,...,v_{g(v)},v_{g(v)+1}\}$ ($0\leq g(v)\leq 5$), defined as follows: $N_{D'}^+(v)=\{v_0,t(v),t'(v)\}$ where $v_0 \rightarrow t(v)$ and $v_0 \rightarrow t'(v)$, $v_1=v$; Set $T(v)=\{t(v),t'(v)\}$. If $d_{D'}^-(v_0)\geq3$, put $g(v)=0$; if not, let $v_2$ be the unique vertex of $D'$ distinct from $v_1$ such that $v_2 \rightarrow v_0$. We have $v_2 \rightarrow v_1$. Again, if $d_{D'}^-(v_1)\geq3$, put $g(v)=1$; otherwise, let $v_3$ be the unique vertex of $D'$ distinct from $v_2$ such that $v_3 \rightarrow v_1$; such a vertex exists by the above lemma. We have $v_2 \rightarrow v_1$, since otherwise we would have a path $p_4$ in $D$.

\vspace*{4mm}

\noindent We may continue this process until meeting a vertex of in-degree at least three in $D'$; call this vertex $v_{g(v)}$, where $g(v)$ is the number of iterations required. Such a vertex exists and $g(v)\leq5$. In fact suppose that $v_1,...,v_5$ are defined as above and $d_{D'}^-(v_i)=2, \forall i\in\{1,2,3,4\}$. By \emph{lemma 2.11} we have $d_{D'}^+(v_i)=2, \forall i\in\{2,3,4,5\}$. If $d_{D'}^-(v_5)=2$, the vertices $v_2,v_3,v_4$ and $v_5$ will be in the same block of $D_m$. The block of $D_m$ containing $\{v_2,v_3,v_4,v_5\}$ cannot be an odd cycle nor complete since $v_2v_5 \notin E(G(D))$ which contradicts \emph{theorem 2.14}. Set $O(v)=\{z\in D' / z\neq v_{g(v)+1}$ and $z \rightarrow v_{g(v)}\}$; we have $z \rightarrow v_{g(v)+1}$ for every $z$ in $O(v)$.

\vspace*{6mm}

\noindent {\bf Lemma 2.16}: \emph{If u and v are two distinct vertices of $D^o$ then $S(u)\cap S(v)=\phi$}.

\vspace*{2mm}

\noindent {\bf Proof}: Let $S(v)=\{t(v),t'(v),v_0,...,v_{g(v)},v_{g(v)+1}\}$ and $S(u)=\{t(u),t'(u),u_0,...,u_{g(u)},$ $u_{g(u)+1}\}$ and suppose that $\exists w \in S(u) \cap S(v)$. If $w=v$ then since $u$ and $v$ are not adjacent we should have $w\notin\{u,u_0,t(u),t'(u)\}$, otherwise we would have $u=v$ or $u\rightarrow v$. So $w=v\in S(u)\setminus\{u,u_0,t(u),t'(u)\}$ so $w=v=u_i$ with $i\geq2$ which is a contradiction since $d_{D'}^+(v)=3$ and $d_{D'}^+(u_i)=2$. So we conclude that $w\neq v$ and similarly $w\neq u$. If $w \in \{v_0,t(v),t'(v)\}$ (i.e. $v\rightarrow w$), then $w \notin T(u)$ since otherwise we would have $v\in N^-(w)=\{u_0,u_1\}\subset S(u)$ (which is a contradiction), and similarly if $w=u_0$ we will have $w\notin T(v)$ and so $u_0=w=v_0$ which implies that $N^+_{D'}(w)=N^+_{D'}(u_0)=T(u)=N^+_{D'}(v_0)=T(v)$ which is also a contradiction. We conclude that if $w \in \{v_0,t(v),t'(v)\}$ then $g(u)\geq2$ and $w=u_i$ with $i\geq2$; and more precisely we have $i=g(u)$ or $i=g(u)+1$, because otherwise we would have $v\in N^-_{D'}(w)=N^-_{D'}(u_i)=\{u_{i+1},u_{i+2}\}\subset S(u)$ which is a contradiction. Since $v\in N^-_{D'}(w)$ and $N^-_{D'}(u_{g(u)+1})=N^-_{D'}(g(u))\setminus\{u_{g(u)+1}\}=O(u)$, $d^+_{D'}(u_{g(u)+1})=2$ and $d^+_{D'}(v)=3$ we conclude that $v\neq u_{g(u)+1}$ and then $v\in O(u)$. Since $v\in O(u)$ then $v\rightarrow u_{g(u)+1}$ and $v\rightarrow u_{g(u)}$, but $u_{g(u)+1}\rightarrow u_{g(u)}$, we can easily conclude that $v_0=u_{g(u)+1}$ and then $N_{D'}^+(v_0)=\{t(v),t'(v)\}=N_{D'}^+(u_{g(u)+1})=\{u_{g(u)},u_{g(u)-1}\}$ which is a contradiction since $u_{g(u)} \rightarrow u_{g(u)-1}$ and $t(v)t'(v) \notin E(G(D))$. So $w \notin \{v_0,v_1,t(v),t'(v)\}$ and similarly $w \notin \{u_0,u_1,t(u),t'(u)\}$, so $\exists i \geq j \geq 2$ such that $w=u_i=v_j$; we will prove by induction on $0\leq l\leq j-2$ that $u_{i-l}=v_{j-l}$: it is true for l=0, suppose that it is true for $l<j-2$ so $u_{i-l}=v_{j-l}$ which implies that $N_{D'}^+(u_{i-l})=\{u_{i-l-1},u_{i-l-2}\}=N_{D'}^+(v_{j-l})=\{v_{j-l-1},v_{j-l-2}\}$, but $u_{i-l-1} \rightarrow u_{i-l-2}$ and $v_{j-l-1} \rightarrow v_{j-l-2}$ so $u_{i-(l+1)}=v_{j-(l+1)}$. Set $l=j-2$, we conclude that $v_2=u_{i-j+2}$ but this implies that $v\in N_{D'}^+(v_2)=N_{D'}^+(u_{i-j+2})\subset S(u)$ which is a contradiction. \hfill $\boxempty$

\noindent {\bf Lemma 2.17}: \emph{Let $L=\{v_{g(v)}/v\in D^o\}$. We have $\forall v\in L, d_{D'}^-(v)=3$ and $\forall v\in D'\setminus L, d_{D'}^-(v)=2$}.

\vspace*{2mm}

\noindent {\bf Proof}: Let $s=|L|=|D^o|$ and $p=|D'\setminus L|=|D'\setminus D^o|$, we have: $$3s+2p\leq \sum_{v \in D^o} d_{D'}^+(v) + \sum_{v \in D'\setminus D^o} d_{D'}^+(v) = \sum_{v \in D'} d_{D'}^+(v)$$ $$|E(D')| = \sum_{v \in D'} d_{D'}^+(v) = \sum_{v \in D'} d_{D'}^-(v)$$ $$\sum_{v \in D'} d_{D'}^-(v) = \sum_{v \in L} d_{D'}^-(v) + \sum_{v \in D'\setminus L} d_{D'}^-(v) \leq 3s+2p$$ So we conclude that all the inequalities are in fact equalities, which holds only if we have $\forall v\in L\;d_{D'}^-(v)=3$, $\forall v\in D'\setminus L\;d_{D'}^-(v)=2$, $\forall v\in D^o\;d_{D'}^+(v)=3$ and $\forall v\in D'\setminus D^o\;d_{D'}^+(v)=2$. \hfill $\boxempty$

\vspace*{6mm}

\noindent {\bf Corollary 2.18}: \emph{For all $v \in D^o$, $O(v)$ contains exactly two vertices}.

\vspace*{2mm}

\noindent {\bf Proof}: Clear, by the definition of $O(v)$, and by \emph{lemma 2.17}. \hfill $\boxempty$

\vspace*{6mm}

\noindent {\bf Proof of Theorem 2.1}: Define the sets: $$S=\bigcup_{v\in D^o}S(v), O=\bigcup_{v\in D^o}O(v), T=\bigcup_{v\in D^o}T(v)$$ We have $|O|\leq|T|$. Suppose that $O=T$, then $D'=D'[S]$ because otherwise we can find a vertex $w$ outside $S$ which is adjacent to a vertex $v$ of $S$ ($D'$ is connected since it is 5-critical) and so $w\in N_{D'}(v)$ which means that $w\in S$ or $w\in O$ (see the definitions of $S(v)$, $O(v)$ and $T(v)$), then since $O=T\subset S$ then in all cases we have $w\in S$ which is a contradiction. Let $v$ be a vertex of $D^o$, then put $c(t(v))=c(t'(v))=1$, $c(v_0)=2$ and $c(v_1)=3$. If $g(v)=0$ we are done, otherwise the colors 1,2 and 3 suffice to color the vertices of $S(v)\setminus\{v_{g(v)},v_{g(v)+1}\}$, let $i\in\{2,3\}\setminus\{c(v_{g(v)-1})\}$ then put $c(v_{g(v)})=4$ and $c(v_{g(v)+1})=i$. We can easily check that $c$ is a good 4-coloring of the 5-chromatic digraph $D'$ which is a contradiction.

\vspace*{4mm}

\noindent So $O\neq T$ which means that $O\nsubseteq T$ or $T\nsubseteq O$. Since $|O|\leq|T|$ then $T\nsubseteq O$, and so we can find a vertex in $T$ which is not in $O$. So we can find a vertex $v\in D^o$ such that $t(v)\notin O(v)$ or $t'(v)\notin O(v)$. We can assume without loss of generality that $t(v)\notin O(v)$ which implies that $N_{D'}^+(t(v))\cap S=\phi$. Let $N_{D'}^+(t(v))=\{u,u'\}$, $\{u,u'\} \cap (D^o\cup L) = \phi$ so $d_{D'}^+(u)=d_{D'}^-(u)=d_{D'}^+(u')=d_{D'}^-(u')=2$. If $u$ and $u'$ are not adjacent, we can find $w \in D'\setminus\{t(v),u,u'\}$ such that $w \rightarrow u$, and if they are adjacent we can assume without loss of generality that $u \rightarrow u'$ and so again we can find $w \in D'\setminus\{t(v),u,u'\}$ such that $w \rightarrow u$. So without losing generality we can say that in all cases we can find $w \in D'\setminus\{t(v),u,u'\}$ such that $w \rightarrow u$. $w \nrightarrow t(v)$ since otherwise $w\in N^-_{D'}(t(v))=\{v_0,v_1\}$ and $u\in N^+_{D'}(w) \subset \{v_0,t(v),t'(v)\}$, but $t(v)\rightarrow u$ and $v_0\rightarrow t(v)$ so $u\neq v_0$ and $u\neq t(v)$, then $u=t'(v)$ which is contradiction since $t(v)t'(v)\notin E(G(D))$, and so $w\nrightarrow t(v)$. Since $d_{D'}^+(w)\geq2$, $\exists w'\in D'\setminus\{u,w,t(v)\}$ such that $w\rightarrow w'$. If $w'\neq u'$, $w'wut(v)u'$ would be a $p_4$ which is a contradiction. We conclude that $N_{D'}^+(w)=\{u,u'\}$, and we have also $w\notin L$ because otherwise we would have $u\in S$. Then $u$,$u'$,$t(v)$ and $w$ are of degree 4, and so they are in the same block of $D_m$ which cannot be neither an odd cycle nor complete which contradicts \emph{theorem 2.14}. \hfill $\boxempty$

\section{Our new shorter proof}

\noindent We provide a new shorter proof of El-Sahili's theorem, which is elementary in the sense that it does not use Gallai's theorem. We will use all the theorems, lemmas and corollaries of the first step.

\vspace*{6mm}

\noindent {\bf New proof of theorem 2.1}: For all $v$ in $D^o$ we define $v'$, $t(v)$ and $t'(v)$, such that $N_{D'}^+(v)=\{v',t(v),t'(v)\}$, $v' \rightarrow t(v)$ and $v' \rightarrow t'(v)$. Let $S(v)=\{v\}\cup N_{D'}^+(v)$, $H(v)=\{v,v'\}$, $O(v)=N_{D'}^-(v')\setminus\{v\}$ and $P(v)=N_{D'}^-(v)\setminus O(v)$. Note that $O(v)$ is not empty since $d_{D'}^-(v')\geq2$ while $P(v)$ can be empty. $\forall w \in O(v), w\rightarrow v$ because otherwise $w'wv'vw''$ would be a $p_4$ where $w'\in N^+_{D'}(w)\setminus\{v',v\}$ and $w''\in \{t(v),t'(v)\}\setminus\{w'\}$. By \emph{lemma 2.12}, every vertex in $O(v)$ has only two out-neighbors i.e. $v$ and $v'$, in particular $O(v)$ is stable.

\vspace*{4mm}

\noindent If $P(v)$ is not empty then $\forall w \in P(v), \exists w' \in O(v)$ such that $w\rightarrow w'$, since otherwise $d^+(w')\geq2$ implies that $\exists w' \in D'\setminus (O(v)\cup\{v,v',w\})$ such that $w\rightarrow w'$ which means that $w'wvuv'$ is a $p_4$ where $u\in O(v)$. By \emph{lemma 2.12}, every vertex in $P(v)$ has only two out-neighbors i.e. $v$ and one vertex in $O(v)$, in particular $P(v)$ is stable.


\vspace*{4mm}

\noindent Let $D^o=\{v_1,v_2,...,v_l\}$, we define $D_i$, $S_i(v)$, $O_i(v)$ and $P_i(v)$ for $0\leq i\leq l$ and $v\in D^o$ as follows: $D_0=D'$, $S_0(v)=S(v)$, $O_0(v)=O(v)$ and $P_0(v)=P(v)$. $D_{i+1}$, $S_{i+1}(v)$, $O_{i+1}(v)$ and $P_{i+1}(v)$ are obtained from $D_i$, $S_i(v)$, $O_i(v)$ and $P_i(v)$ by removing $S_i(v_{i+1})$ and then contracting $O_i(v_{i+1})$ and $P_i(v_{i+1})$ if any of them is not empty.

\vspace*{4mm}

\noindent We can easily check that all the vertices of $D_l$ has at most two out-neighbors. Suppose that $D_l$ contains a 5-tournament $T$, then $T$ contains at least one contracted vertex $w$ (Otherwise $T$ would be a subdigraph of $D'$). $w=v_{O_{l-1}(v)}$ or $v_{P_{l-1}(v)}$ for some $v\in D^o$, and in both cases $w$ has at most one out-neighbor in $D_l$, and this means that: $$10=|E(T)|=\sum_{u\in T} d_{T}^+(u)=d_{T}^+(w) + \sum_{u\in T\setminus\{w\}} d_{T}^+(u)\leq 1 + 4 \times 2=9$$ Which gives a contradiction. So $D_l$ does not contain any 5-tournament and by \emph{corollary 2.5} we conclude that $\chi(D_l)\leq4$.

\vspace*{4mm}

\noindent Let $i$ be the least integer such that $\chi(D_i)\leq4$, then $i>0$ because $D_0=D'$ is 5-critical. Let $v=v_i$ and let $c$ be a good 4-coloring of $D_i$. Color the vertices in $O_{i-1}(v)$ by $c(v_{O_{i-1}(v)})$ and color those in $P_{i-1}(v)$ by $c(v_{P_{i-1}(v)})$.

\vspace*{4mm}

\noindent If $t(v)$ and $t'(v)$ are adjacent (in $D_{i-1}$) to at most three colors, we color them by a remainder color, then similarly color $v$ and then $v'$ (They are each adjacent to at most three colors) and we get $\chi(D_{i-1})\leq4$, a contradiction.

\vspace*{4mm}

\noindent We conclude that $t(v)$ and $t'(v)$ are adjacent to the four colors 1,2,3 and 4. We may assume without loss of generality that $t(v)$ is adjacent to 1 and 2 and that $t'(v)$ is adjacent to 3 and 4. If $O_{i-1}(v)=\phi$, color $t(v)$ by 3, $t'(v)$ by 1, $v$ by 2 and $v'$ by 4, and we get good 4-coloring. So $O_{i-1}(v)\neq\phi$, we may assume without loss of generality that $c(v_{O(v)})=1$. Color $t(v)$ by 3 and color $t'(v)$ by 1, $v$ by 2 and $v'$ by 4, and we get good 4-coloration and $\chi(D_{i-1})\leq4$, a contradiction. \hfill $\boxempty$

\section{Conclusion}

\noindent In this chapter we have presented El-Sahili's theorem ~\cite{pa1} stating that we can always find a copy of the anti-directed path $p_4$, in any $5$-chromatic digraph where every vertex has at least two out-neighbors and which is not exactly $T_5$. We have presented El-Sahili's proof and we have provided a new shorter proof.

\vspace*{4mm}

\noindent Is the condition that every vertex has at least two out-neighbors really necessary? El-Sahili gave a positive answer in his paper through the following example: Construct a digraph by adding to $T_5$ an arc $(x,y)$ where $x\notin T_5$ and $y\in T_5$, then we can easily check that this digraph does not contain a copy of $p_4$.

\vspace*{4mm}

\noindent The example given above contains $T_5$ and this shows that the condition that every vertex has at least two out-neighbors is necessary for digraphs containing $T_5$. What if it does not contain $T_5$? El-Sahili concluded his paper ~\cite{pa1} by asking the following question: \emph{Can we find a 5-chromatic digraph which contains neither a 5-tournament nor $p_4$?}

\vspace*{2mm}

\noindent We conclude this chapter by stating the following conjecture of us:

\vspace*{6mm}

\noindent {\bf Conjecture 2.19}: \emph{Let $D$ be a $2n+1$-chromatic graph where $n\geq2$. If $D$ does not contain any $2n+1$-tournament, and if every vertex of $D$ has at least $n$ out-neighbors. Then $D$ contains the antidirected path $p_{2n}$ of length $2n$ starting with a backward arc}.

\chapter{Paths with two blocks in $n$-chromatic digraphs}

\section{Introduction}

\noindent An important problem in graph theory is to find which oriented paths can be found in $n$-chromatic digraphs. Gallai-Roy's celebrated theorem ~\cite{pa8, pa9} states that every $n$-chromatic digraphs contains a directed path of length $n-1$. The question is that can we find an oriented path of length $n-1$ with more than one block? or more generally, how big should be the chromatic number of a digraph to guarantee the existence of an oriented path of length $n-1$?

\vspace*{4mm}

\noindent Burr ~\cite{pa11} proved that every $(n-1)^2$-chromatic digraph contains any tree of order $n$, in particular every $(n-1)^2$-chromatic digraph contains any oriented path of length $n$. In this chapter we are interested in paths with two blocks. El-Sahili ~\cite{pa15} introduces the function $f(n)$ which is defined to be the minimal integer $f(n)$ such that every $f(n)$-chromatic digraph contains any path with two block $P(k,l)$ with $k+l=n-1$, and he conjectured that $f(n)=n$ for $n\geq 4$. El-Sahili proved ~\cite{pa15} that $f(n)\leq\frac{3}{4}n^2$. El-Sahili and Bondy ~\cite{pa15} proved that the conjecture holds when one of the two blocks have length 1.

\vspace*{4mm}

\noindent El-Sahili and Kouider ~\cite{pa16} introduced the notion of maximal spanning out-forest and used it to prove that $f(n)\leq n+1$. Addario-Berry et al ~\cite{pa2} used strongly connected digraphs and maximal spanning out-forests to prove El-Sahili's conjecture ($f(n)=n$ for $n\geq4$). Later El-Sahili and Kouider ~\cite{pa3} provided a new elementary proof of El-Sahili's conjecture without using strongly connected digraphs. In this chapter we provide a detailed explanation of both methods. We show that the first method contains a small error and we provide a correction.

\section{Maximal spanning out-forest}

\noindent The level $l_F(v)$ of a vertex $v$ in an out-forest $F$ is defined as in the case of out-branching; the order of a longest directed path ending at $v$. We denote by $T_v(F)$ the out-branching of $F$ rooted at $v$ and by $P_v$ the directed path in $F$ of order $l_F(v)$ which ends at $v$. For all $u \in P_v$, $P_uv$ denotes the $uv$-directed path in $F$.

\vspace*{4mm}

\noindent Let $D$ be a digraph, a spanning subdigraph $F$ of $D$ is said to be a \emph{maximal spanning out-forest} if $F$ is a out-forest such that $\forall x,y \in V(D)$, if $x\rightarrow y$ with $l_F(x)\geq l_F(y)$ then there exists a directed path from $y$ to $x$ in $F$, i.e. $y\in P_x$. The set $L_i$ of vertices having the same level $i$ is a stable (by definition).

\vspace*{4mm}

\noindent Let $F$ be an out-forest which is a spanning subdigraph of a digraph $D$. If $F$ is not a maximal out-forest, then there exist an arc $x\rightarrow y$ such that $l_F(x)\geq l_F(y)$ and there is no directed path from $y$ to $x$ in $F$, the out-forest $F'$ obtained from $F$ by deleting the arc whose head is $y$ (If such one exists) and adding the arc $x\rightarrow y$ is called an \emph{elementary improvement} of $F$.

\vspace*{4mm}

\noindent We can easily see that the level of each vertex in $F'$ is at least its level in $F$, and there exists a vertex ($y$) whose level strictly increases. Since the level of a vertex cannot increase infinitely (The maximum level that can be reached is |V(D)|), we can see that after a finite number of elementary improvements we get to a maximal spanning out-forest which is call a \emph{maximal closure} of $F$. Thus starting with a spanning out-forest that contains no arcs we can prove the existence of a maximal spanning out-forest of $D$. We have also another way to get the existence of a maximal spanning out-forest; choose an out-forest $F$ which maximizes the sum of the levels of all vertices.

\vspace*{4mm}

\noindent The notion of maximal spanning out-forests introduced by El-Sahili and Kouider ~\cite{pa16} is useful in the context of universal digraphs. As shown by El-Sahili and Kouider ~\cite{pa16}, it gives an easy proof of Gallai-Roy's theorem. Indeed, consider a maximal spanning out-forest of an $n$-chromatic digraph $D$. Since every level is a stable set, there are at least $n$ levels. Hence $D$ contains a directed path of length at least $n-1$. Final forests are also useful for finding paths with two blocks, as illustrated by the following proof due to El-Sahili and Kouider ~\cite{pa16}.

\vspace*{6mm}

\noindent {\bf Lemma 3.1} ~\cite{pa16}: \emph{Let $F$ be a maximal spanning out-forest of a digraph $D$. If $v\rightarrow w$ is an arc from $F_i$ to $F_j$. Then}
\begin{enumerate}
  \item \emph{If $k \leq i < j - l$, then $D$ contains a $P(k, l)$}.
  \item \emph{If $k < j \leq i - l$, then $D$ contains a $P(k, l)$}.
\end{enumerate}

\vspace*{2mm}

\noindent {\bf Proof}: 1. Let $P_l$ be the directed path of $F$ which starts at $F_{j-l}$ and ends at $w$ and $P_{k-1}$ be the directed path in $F$ starting at $F_{i-(k-1)}$ and ending at $v$. Then $P_{k-1} \cup vw \cup P_l$ is a $P(k, l)$.

\vspace*{4mm}

\noindent 2. Let $P_{l-1}$ be the directed path in $F$ which starts at $F_{i-l+1}$ and ends at $v$. Let $P_k$ be the directed path in $F$ starting at $F_{j-k}$ and ending at $w$. Then $P_k \cup P_{l-1} \cup vw$ is a $P(k, l)$. \hfill $\boxempty$

\vspace*{6mm}

\noindent {\bf Corollary 3.2} ~\cite{pa16}: \emph{Every digraph with chromatic number at least $k+l+2$ contains a $P(k, l)$.}

\vspace*{2mm}

\noindent {\bf Proof}: 1. Let $F$ be a maximal spanning out-forest of $D$. Color the levels $F_1,...,F_k$ of $F$ with colors $1,...,k$. Then color the level $F_i$, where $i > k$, with color $j\in\{k+1,...,k+l+1\}$ such that $j \equiv i$ mod $l + 1$. Since this is a $k+l+1$-coloring, it's not a good, and so there exists an arc which satisfies the hypothesis of \emph{Lemma 2.3}. \hfill $\boxempty$

\section{Paths with two blocks in strongly connected digraphs}

\noindent {\bf Theorem 3.3} ~\cite{pa17}: \emph{Every strongly connected digraph $D$ has a circuit of length at least $\chi(D)$.}

\vspace*{2mm}

\noindent Let $k$ be a positive integer and $D$ be a digraph. A directed circuit $C$ of $D$ is $k$-good if $|C| \geq k$ and
$\chi(D[V(C)]) \leq k$. Note that Theorem 3.3 states that every strongly connected digraph D has a $\chi(D)$-good
circuit.

\vspace*{6mm}

\noindent Note that the last part of the proof in ~\cite{pa2} of the following lemma contains an error. We will show the proof in ~\cite{pa2} and explain why it is false, and then we will provide a correction.

\vspace*{6mm}

\noindent {\bf Lemma 3.4}: \emph{Let $D$ be a strongly connected oriented multi-graph and $k\in\{3,...,\chi(D)\}$. Then $D$ has a $k$-good circuit.}

\vspace*{2mm}

\noindent {\bf Proof}: By Bondy's theorem, there exists a circuit with length at least $\chi(D)$, so the lemma is true for $k = \chi(D)$. If $k=3$ then if $C$ is the shortest circuit of $D$, then it's chordless and therefore $\chi(C)=2$ or $3$. Suppose that $3 \leq k < \chi(D)$ and consider a shortest circuit $C$ with length at least $k$. We claim that $\chi(D[V(C)]) \leq k$.Suppose to the contrary that $\chi(D[V(C)]) \geq k + 1$, and let $D'$ be a maximal sub-oriented-graph of $D[V(C)]$ such that $D'$ is a strongly connected digraph in which $C$ is a subdigraph. If any two vertices of $D'$ are adjacent in $D$, they are still adjacent in $D'$, and so $\chi(D')=\chi(D)\geq k+1$, moreover $C$ is a hamiltonian circuit of $D'$.

\vspace*{4mm}

\noindent Let $u$ be a vertex of $D'$, if $v_1,...,v_{k-1}$ are in-neighbors of $u$ in $D'$, listed in such a way that $v_1,...,v_{k-1},u$ appear in the same order along $C$, the sub-circuit of $C+v_{k-2}u$ not containing $v_{k-1}$ would have length at least $k$ since it contains $v_1,...,v_{k-2}$ and $u$ in addition to the out-neighbor of $u$ in $C$. This contradicts the minimality of $C$, so we conclude that every vertex has at most $k-2$ in-neighbors in $D'$ and similarly at most $k-2$ out-neighbors in $D'$.

\vspace*{4mm}

\noindent A \emph{handle decomposition} of $D'$ is a sequence $H_1,...,H_r$ such that:

\vspace*{2mm}

\begin{enumerate}
  \item $H_1$ is a circuit of $D'$.
  \item For $2\leq i\leq r$, $H_i$ is a \emph{handle}, that is, a directed path in $D'$ (with possibly the same end-vertices i.e. a circuit) meeting $V(H_1\cup...\cup H_{i-1})$ exactly at his end-vertices.
  \item $D' = H_1\cup...\cup H_r$.
\end{enumerate}

\vspace*{4mm}

\noindent An $H_i$ which is an arc is a \emph{trivial handle}. It is well-known that $r$ is invariant for all handle decompositions of $D'$ (indeed, $r$ is the number of arcs minus the number of vertices plus one, it is proved by a simple induction on $r$). However the number of nontrivial handles is not invariant. Let us then consider $H_1,...,H_r$, a handle decomposition of $D'$ with minimum number of trivial handles. Since the trivial handles does not add any new vertices, we can enumerate first the nontrivial handles, and so we can assume that $H_1,...,H_p$ are not trivial and that $H_{p+1},...,H_r$ are arcs.

\vspace*{4mm}

\noindent Let $D^o := H_1\cup...\cup H_p$. Clearly $D^o$ is a strongly connected spanning subdigraph of $D$. Observe that since $\chi(D')>3$, $D'$ is not an induced circuit which means that $r>1$, so $p>1$ because otherwise a trivial handle would be a chord of $H_1$ so by shortcutting $H_1$ through this chord we get two non trivial non handles which contradicts the maximality of $p$.

\vspace*{4mm}

\noindent We denote by $x_1,...,x_q$ the handle $H_p$ minus its end-vertices.

\vspace*{4mm}

\noindent If $q=1$, the digraph $D^o-x_1$ is strongly connected, and therefore $D'-x_1$ is also strongly connected. Moreover since $\chi(D')\geq k+1$ we have $\chi(D'-x_1)\geq k$. Thus by Bondy's theorem, there exists a circuit of length at least $k$ in $D'-x_1$ that is shorter than $C$, a contradiction with the minimality of $C$.

\vspace*{4mm}

\noindent If $q=2$, $x_2$ is the unique out-neighbor of $x_1$ in $D'$ because otherwise we would make two non trivial handles out of $H_p$, contradicting the maximality of $p$. Similarly, $x_1$ is the unique in-neighbor of $x_2$. Since the out-degree and the in-degree of every vertex is at most $k-2$, both $x_1$ and $x_2$ have degree at most $k-1$ in the underlying graph of $D$. Since $\chi(D)>k$, it follows that $\chi(D-\{x_1,x_2\})>k$ because otherwise we can extend a good $k$-coloring of $D-\{x_1,x_2\}$ by giving each of $x_1$ (we can always find such a color since $x_1$ is adjacent to at most $k-1$ vertices) and then we do the same with $x_2$. Since $D-\{x1, x2\}$ is strongly connected, it contains, by Bondy's theorem, a circuit with length at least $k$, contradicting the minimality of C.

\vspace*{4mm}

\noindent Hence, we may assume that $q>2$. $\forall i\in\{1,...,q-1\}$, by the maximality of $p$, the unique arc in $D'$ leaving $\{x_1, . . . , x_i\}$ is $x_ix_{i+1}$ (otherwise we would make two nontrivial handles out of $H_p$). Similarly, $\forall i\in\{2,...,q\}$, the unique arc in $D'$ entering $\{x_j,...,x_q\}$ is $x_{j-1}x_j$. In particular, as for $q = 2$, $x_1$ has out-degree 1 in $D'$ and $x_q$ has in-degree 1 in $D'$.

\vspace*{4mm}

\noindent The next paragraph is, word by word (with exception for the terminology), the last part of the proof in ~\cite{pa2} \textbf{which contains an error}:

\vspace*{4mm}

\noindent \emph{``Another consequence is that the underlying graph of $D'-\{x_1, x_q\}$ has two connected components $D_1 = D'-\{x_1,x_2,...,x_q\}$ and $D_2 = D'[\{x_2,...,x_{q-1}\}]$. Since the degrees of $x_1$ and $x_q$ in the underlying graph of $D'$ are at most $k-1$ and $D'$ is at least $(k+1)$-chromatic, it follows that $\chi(D_1)$ or $\chi(D_2)$ is at least $k + 1$. Each vertex has in-degree at most $k-2$ in $D'$ and $d_{D_2}^+(x_i) \leq 1$ for $2\leq i\leq q-1$, so $\Delta(D_2) \leq k-1$ and $\chi(D_2)\leq k$. Hence $D_1$ is at least $(k + 1)$-chromatic and strongly connected. Thus by Bondy's theorem, $D_1$ contains a circuit of length at least $k$ but shorter than $C$. This is a contradiction."} ~\cite{pa2}

\vspace*{4mm}

\noindent \textbf{The error} is that there is no reason to say that $d_{D_2}^+(x_i) \leq 1$ for $2\leq i\leq q-1$, in fact $x_i \nrightarrow x_j$ for $j>i$ but we can have $x_i \rightarrow x_j$ for $j<i$ and so we can have $d_{D_2}^+(x_i) > 1$, and therefore $\Delta(D_2)$ can be greater than $k-1$. So we will prove that $\chi(D_2)\leq k$ through another proof:

\vspace*{4mm}

\noindent Let $D^i := D[\{x_i,...,x_{q-1}\}]$ and let $i$ be the minimum integer greater than 1 such that $\chi(D^i) \leq k$. Suppose that $i>2$: since the unique arc in $D$ entering $\{x_i , ... , x_q\}$ is $x_{i-1}x_i$ then we have $d^+_{D^{i-1}}(x_{i-1}) = 1$ and since $d^-_D(x_{i-1}) \leq k-2$ we have $d_{D^{i-1}}(x_{i-1}) \leq k-1$ and therefore $\chi(D^{i-1}) \leq k$ which contradicts the minimality of $i$. Then $i=2$ and $\chi(D_2) = \chi(D^2) \leq k.$ \hfill $\boxempty$

\vspace*{6mm}

\noindent The existence of good circuits directly implies the main theorem in the case of strongly connected digraphs.

\vspace*{6mm}

\noindent {\bf Corollary 3.5}: \emph{Let $k + l = n - 1$ where $n\geq4$ and let $D$ be a strongly connected $n$-chromatic digraph then $D$ contains a $P(k,l)$.}

\vspace*{2mm}

\noindent {\bf Proof}: Since $P(k,l)$ and $P(l,k)$ represent the same digraph and since $k+l=n-1\geq3$, we may assume that $l\geq(n-1)/2\geq3/2$ which means that $l\geq2$. By \emph{lemma 3.4} $D$ contains an $(l+1)$-good circuit $C$, the chromatic number of the (strongly connected) contracted oriented multi-graph $D/C$ is at least $k$, since otherwise we may use a good $k$-coloring of $D/C$ to construct a good $n-1$-coloring of $D$: keep the colors of the vertices of $D-C$, and for the vertices of $C$ we give one vertex the color of $v_C$ and then we color the other vertices by $l$ new colors. We conclude that $\chi(D/C)\geq k+1$ and by Bondy's theorem, $D/C$ has a circuit of length at least $k+1$, and in particular the vertex $v_C$ is the end of a path $P$ of length $k$ in $D/C$. Finally $P\cup C$ contains a $P(k, l)$. \hfill $\boxempty$

\section{General case, first method (Addario-Berry et al)}

\noindent {\bf Theorem 3.6}: \emph{Let $k+l=n-1\geq3$ and let $D$ be an $n$-chromatic digraph. Then $D$ contains a $P(k,l)$.}

\vspace*{2mm}

\noindent {\bf Proof}. We can assume that $l\geq k$, and therefore $l\geq2$. Suppose to the contrary that $D$ does not contain $P(k,l)$. Let $F$ be a maximal spanning out-forest of $D$.

\noindent Consider the following coloring (Which we call \emph{canonical coloring}) of $D$: for $1\leq i \leq k-1$, the vertices of $F_i$ are colored $i$, and for $i\geq k$, the vertices of $F_i$ are colored $j$, where $j\in\{k,...,k+l\}$ and $j\equiv i$ mod $l + 1$. Since we colored $D$ with less than $n$ colors, this coloring can not be good. In particular, there exists an arc $v\rightarrow w$ from $F_i$ to $F_j$ where $i,j \geq k$ and $j \equiv i$ mod $l+1$. By \emph{Lemma 3.1 (1)}, we get a contradiction if $i<j$. Thus $j < i$, and by \emph{Lemma 3.1 (2)}, we necessarily have $j = k$ and $i \geq k + l + 1$. Since $F$ is a maximal spanning out-forest we can find in $F$ a directed path from $w$ to $v$. In particular $F+vw$ has a circuit $C$ of length at least $l+1$. If $\chi(D[C])\leq l+1$ then $C$ is $(l+1)$-good, if not, then by \emph{Lemma 3.4}, it contains an $(l+1)$-good circuit. So in all cases we can find an $(l+1)$-good circuit which is disjoint from $F_1\cup...\cup F_{k-1}$.

\vspace*{4mm}

\noindent We inductively define couples $(D^i, F^i)$ as follows: Set $D^0:=D$, $F^0:=F$. Then, if there exists an $(l+1)$-good circuit $C^i$ of $D^i-F_1^i\cup...\cup F_{k-1}^i$, define $D^{i+1}:=D^i-V(C^i)$ and let $F^{i+1}$ be any maximal closure in $D^{i+1}$ of $F_i-V(C_i)$.

\vspace*{4mm}

\noindent With the previous definitions, we have $D^1 = D - V (C^0)$. This inductive definition certainly stops on some $(D^p, F^p)$ where the canonical coloring of $D^p$ is a good coloring.

\vspace*{4mm}

\noindent At each inductive step, the circuit $C^i$ must contain a vertex $v^i$ of $F^i_k$, otherwise the union of $C^i$ (which has length at least $l+1$) and a path of $F^i$ starting at $F^i_1$ and ending at $C^i$ (which would have length at least $k$ if $C^i$ does not meet $F^i_k$) would certainly contain a $P(k,l)$. Let $u^i$ the unique in-neighbor of $v^i$ in $F^i_{k-1}$. $\forall j>i,\;l_{F^j}(u^i)=k-1$, since $l_{F^j}(u^i)\geq k-1$ because we apply successive elementary improvements, and $l_{F^j}(u^i)$ cannot be greater than $k-1$, otherwise $u^i$ would be the end of a path $P$ of length $k-1$ in $D-C^i$ and thus $Ci\cup P\cup u^iv^i$ would contain a $P(k,l)$. Thus every circuit $C^i$, $i=0,...,p-1$, has an in-neighbor $u^i$ in $F^p_{k-1}$.

\vspace*{4mm}

\noindent Observe that we cannot have any arc between two circuits $C^i$ since they are disjoint and the length of each one is at least $l+1$, and if there is such an arc we get a $P(k,l)$ since $l\geq k$. Observe also that no vertex of $C^i$ has a neighbor, (in- or out-), in any level $F^p_j$ for any $j>k$ because otherwise we get a $P(k,l)$. Moreover, no vertex of $C^i$ has an in-neighbor in $F^p_k$.

\vspace*{4mm}

\noindent Let us call \emph{bad vertices} the out-neighbors of the vertices of all $C^i$ in $F^p_k$ and \emph{good vertices} the non-bad vertices in $F^p_k$. A bad vertex $b$ cannot have in-neighbors in more than one circuit $C^i$, since the length of those circuits is at least $l+1$ and so joining two circuits $C^i$ with $b$ through two arcs towards it make a $P(k,l)$. Moreover $b$ has at most $l$ in-neighbors in $C^i$: Suppose to the contrary that $w_1,...,w_{l+1}$ are in-neighbors of $b$ in $C^i$, enumerated with respect to the cyclic order of $C^i$ such that $w_1$ is the first vertex $w_j$ along $C^i$ which appears after $v^i$ (i.e. $C^i[v^i,w_1] \cap \{w_1,...,w_{l+1}\}=\{w_1\}$). Let $P$ be the path of $F^p$ starting at $F^p_1$ and ending at $u^i$. Now $P\cup u^iv^i \cup C[v^i,w_1] \cup w_1b \cup C[w_2,w_{l+1}] \cup w_{l+1}b$ contains a $P(k, l)$, a contradiction.

\vspace*{4mm}

\noindent Let $b$ is a bad vertex, we denote by $S_b$ the set of \emph{descendants} of $b$ in $F_p$, i.e. the set of vertices $x$ such that there is a path from $b$ to $x$ in $F_p$, including $b$ itself.

\vspace*{4mm}

\noindent We claim that every arc $x\rightarrow y$ entering $S_b$ (i.e. $y\in S_b$ and $x \notin Sb$) in $D' := D - F^p_1\cup...\cup F^p_{k-1}$ is such that $y=b$ and $x\in C^i$. Indeed, suppose that $y\neq b$, $y$ would be a strict descendant of $b$ in $F^p$ and then $l_{F^p}(y)> k$ and so $x\notin C^j\; \forall j\in\{1,2,...,p-1\}$, thus $x \in F^p$. Let $P_1$ be the path in $F^p$ (of length at least $k-1$) ending at $x$, let $P_2$ be the path in $F^p$ starting at $b$ and ending at $y$ and let $v$ be an in-neighbor of $b$ in $C^i$. $P_1\cup xy\cup C^i\cup vb \cup P_2$ would contain a $P(k,l)$, which gives a contradiction; so we conclude that $y=b$, if $x\in F^p$ we would have $x\rightarrow b$, $l_{F^p}(x)\geq l_{F^p}(b)=k$ without having any directed path from $b$ to $x$ which contradicts the fact that $F^p$ is a spanning maximal out-forest. So we must have $x\notin F^p$ which means that there exist $0\leq j<p$ such that $x\in C^j$. We must have $j=i$ since $b$ cannot have in-neighbors in more than one circuit $C^j$.

\vspace*{4mm}

\noindent We claim also that we have no arcs leaving $S_b$. Indeed, let $x\rightarrow y$ be an arc of $D'$ such that $x \in S_b$ and $y \notin S_b$. If $y \in F^p$, there exists a path $P_1$ (of length at least $k4$) in $F^p$ ending at $y$ which does not meet $S_b$ nor $C^i$. Let$P_2$ be the path in $F^p$ which starts at $b$ and ends at $x$, and let $v$ be an in-neighbor of $b$ in $C^i$. We can then find a copy of $P(k,l)$ in $P_1 \cup C^i \cup vb \cup P_2 \cup xy$. Thus $y\notin F^p$ and therefore it belongs to some $C^j$ , but this is impossible since $l_{F^p}(y)\geq k$.

\vspace*{4mm}

\noindent We resume:
\begin{itemize}
  \item There is no arcs between different $C^i$'s.
  \item Each $C^i$ is adjacent to a unique vertex in $D'$ which is bad.
  \item If $b$ is a bad vertex, then the only arcs between $S_b$ and $D'-S_b$ are those between $b$ and a unique circuit $C^i$, we have at most $l$ such arcs.
\end{itemize}

\vspace*{4mm}

\noindent Let us color $D$ with $n-1$ colors. Let $D_1$ be the subdigraph $D$ induced by the vertices of $F^p$ which are not in $S_b$ for any bad vertex $b$. The canonical coloring of $D_1$ is good since all the vertices in $D_1$ of level $k$ are good. We will extend this coloring for the other vertices of $D$ (which are vertices of some $C^i$, or descendants of some bad vertex).

\vspace*{4mm}

\noindent Every $C^i$ is $(l + 1)$-good and thus $(l + 1)$-colorable. Moreover, we have no arcs between any two circuits $C^i$, so we may color their vertices by the colors $k,k+1,...,k+l$. This extension of the coloring is also good since the vertices whose level is at most $k-1$ are colored with colors $1,...,k-1$, and the vertices of $D_1$ whose level is at least $k$ are descendants of good vertices.

\vspace*{4mm}

\noindent So it remains to extend the coloring for the descendants of bad vertices. Let $b$ be a bad vertex, then $b$ is adjacent (in $D'$) to at most $l$ vertices in some unique $C^i$, so we can properly choose a color $c$ for $b$ from the $l+1$ colors $k,k+1,...,k+l$. Since the strict descendants of $b$ are not adjacent to any vertex outside $S_b$, we properly color any descendant $v$ of $b$ with a color $c(v)$ in $\{k,k+1,...,k+l\}$ such that $c(v)\equiv c+l_{f^p}(v)$ mod $l+1$. We get a good $n-1$-coloring of $D$, which is a contradiction. \hfill $\boxempty$

\section{General case, second method (El-Sahili and Kouider)}

\noindent To prove \emph{theorem 3.6}, we will use the following weaker result, proved by El-Sahili and Bondy:

\vspace*{6mm}

\noindent {\bf Theorem 3.7} ~\cite{pa15}: \emph{For $n \geq 4$, every $n$-chromatic digraph contains a path $P(n - 2, 1)$.}

\vspace*{6mm}

\noindent We explain now the new method of El-Sahili and Kouider to prove \emph{theorem 3.6}:

\vspace*{2mm}

\noindent {\bf New proof of theorem 3.6}. Let $D$ be an $n$-chromatic digraph. Due to \emph{theorem 3.7}, it is sufficient to prove that $D$ contains any path $P(k, l)$ with $2\leq k \leq l$ and $k + l = n - 1$. Consider a maximal spanning out-forest $F$ of $D$ minimizing $u_{k}(F) = \sum_{j=1}^{k-1} |L_j(F)|$ . The vertices in $U_{i} = L_{i}(F)$ are taken the color $i$ for $1 \leq i \leq k - 1$. For $i \leq l$, set $U_{k+i} = \cup_{r\geq0} L_{k+i+r(l+1)}(F)$.

\vspace*{4mm}

\noindent Step 1: Suppose to the contrary that $D$ contains no path $P(k, l)$. Then $U_i$ is a stable set for $i \neq k$. Indeed, this fact is trivial for $i \leq k - 1$. If $U_{i}$ is not stable for $i > k$, then there is an edge $uv \in G(D[U_{i}])$. Since vertices having the same level are not adjacent, we must have $l_{F} (u)\neq l_{F}(v)$, then $|l_{F}(u) - l_{F}(v)| \geq l + 1$ and $min(l_{F}(u), l_{F}(v))\geq k + 1$, so by \emph{lemma 3.1} $D$ contains a path $P(k, l)$ which is a contradiction. if $U_k$ is stable, we get $n-1$ stables which contradicts the fact that $\chi(D)=n$, then $U_k$ is not stable. By \emph{lemma 3.1} the only possible arcs in $U_{k}$ are those with heads in $L_k(F)$. These vertices of $L_{k}(F)$ are said to be bad. It is clear that if $v$ is a bad vertex then $T_{v}(F)$ contains a circuit of length at least $l + 1$, and so each vertex in $T_{v}(F)$ is the end of a directed path of length $l$, and this means that:

\vspace*{4mm}

\noindent There is no edge $uw$ in $G(D)$ with $u \in T_{v}(F)$ and $w \notin T_{v}(F)$ such that $l_F(w) \geq k$ (1).

\vspace*{4mm}

\noindent We get a contradiction if we give the uncolored vertices colors in ${1,...,k,...,k+l}$ to obtain a good $(n-1)$-coloring of $D$. By remark (1) This can be done separately on each $T_{v}(F)$ where $v$ is bad. Let $v$ be a bad vertex of $F$ and suppose that $F$ is chosen as above with a minimal number of bad vertices.

\vspace*{6mm}

\noindent Step 2: Let $x, y \in N^-(v) \cap U_k$, we have $l_F (x) = l_F (y)$ since otherwise we will have $l_F (x) - l_F (y) \geq l + 1$, and so $P_{vy} \cup P_{y'x} \cup xv \cup yv$, where $y' \in P_x$ and $l_F (y') = l_F (y) + 1$, contains a path $P(k, l)$. Set $h(v) = l_F (x) = l_F (y)$. A vertex $u \in D$ is said to be rich in $F$ if $l_F (u) \geq k$ and $N(u) \cap L_i(F) \neq \phi$ for all $i \leq k-1$. If $N^-(v)\cap U_k$ contains no rich vertices, then each vertex $u \in N^-(v)\cap U_k$ can take a color $i \leq k - 1$ such that $N(u) \cap L_i(F) = \phi$. A remainder vertex $x$ takes the color $k \leq i \leq k + l$ if $x \in U_i$. We obtain an good $(n - 1)$-coloring which is a contradiction. Similarly we verify that $v$ is rich. Let $u$ be a rich vertex in $N^-(v) \cap U_k$. we have $N^-(v) \cap U_k = {u}$. In fact if there is another vertex $w \in N^-(v) \cap U_k$, let $s$ be the smallest integer such that $N^+(u) \cap L_s(F) \neq \phi$, we have $s \leq k$. Let $x \in N^+(u) \cap L_s(F)$. Since $F$ is a maximal spanning out-forest then $x \in P_u$ which contains $P_v$. If $s = 1$, $ux \cup P_v \cup P_{vw} \cup wv$ contains a path $P(k, l)$. If $s > 1$, then let $y \in N^-(v) \cap L_{s-1}(F)$, $y$ exists due to the minimality of $s$, $P_y \cup yu \cup ux \cup P_{xv} \cup P_{vw} \cup wv$ contains a path $P(k, l)$. The same argument proves that:

\vspace*{4mm}

\noindent $u$ is the unique rich in-neighbor of $v$ with level greater than $n - 1$ (2).

\vspace*{4mm}

\noindent Denote by $\overline{v}$ the vertex $u$ and by $C_v$ the circuit $P_{vu} \cup uv$ and set $C_v = v_kv_{k+1}...v_pv_k$ where $v_k = v$ and $v_p = u$. Note that there exist an integer $f$ such that $l(C_v)=1+f(l+1)$. We show that $v_{k+1}$ is a rich vertex: $N(v) \cap U_{k+1}$ must contain a rich vertex $x$, because otherwise we may give all the vertices in $N(v) \cap U_{k+1}$ an appropriate color in $\{1,2,...,k-1\}$ and then give $v$ the color $k+1$, and the color $i$ for remaining vertices in $T_v \cap L_i$. We get then a good $n-1$-coloring, a contradiction. Then we must have $x \in N^+(v) \cap L_{k+1}(F)$ by remark (2). If $v_{k+1}$ is not rich then $x \notin C_v$. We show as above that $N(x) \cap L_i(F) = N^-(x) \cap L_i(F)$ for all $i \leq k - 1$: If $\exists s, N^+(x) \cap L_i(F) \neq \phi$, we may suppose that $s$ is minimal, let $y\in N^+(x) \cap L_s(F)$. If $s = 1$, $xy \cup P_v \cup C_v$ contains a path $P(k, l)$. If $s > 1$, then let $y' \in N^-(x) \cap L_{s-1}(F)$, $y'$ exists due to the minimality of $s$, $P_{y'} \cup y'x \cup xy \cup P_{yv} \cup C_v$ contains a path $P(k, l)$.

\vspace*{4mm}

\noindent If $zw \in E(G(D))$ with $w \in T_v - T_x$ and $z \in T_x$, we have $w = v$ and $z = x$: Suppose to the contrary that $z\neq x$, since $V(C_v) \subset V(T_v - T_x)$ $w$ is the end of a directed path $Q_w$ of length $l$ in $T_v - T_x$. Let $y \in N^-(x) \cap L_{k-1}(F)$ so $P_y \cup yx \cup P_{xz} \cup Qw \cup wz$ contains a path $P(k, l)$ regardless of the orientation of $wz$. So we $z=x$, but this means that $w\rightarrow z$ because otherwise $P_y \cup yx \cup xw \cup Qw$ would again contain a $P(k,l)$. $w\rightarrow z$ means that either $l_F(w)<l_F(z)$ or $z\in T_w$, the latter case does not hold since $w \in T_v - T_x$ and $z \in T_x$, so $l_F(w)<l_F(z)$ and thus $w$ is necessarily $v$. We conclude that $vx$ is the only edge between $T_v-T_x$ and $T_x$. (3)

\vspace*{4mm}

\noindent Color a vertex $z \in T_x \cap U_i$ by the color $i + 1$ if $i < n - 1$ and by the color $k$ if $i = n - 1$. We do the same with any other rich neighbor of $v$ in $U_{k+1}$. We give the other vertices of $N(v) \cap U_{k+1}$ appropriate colors from $\{1,...,k-1\}$, $v$ is then colored by $k + 1$ and each remainder vertex $z \in U_i$ ($k + 1 \leq i \leq k + l$) is colored by the color $i$. We get an good $(n - 1)$-coloration of $D$, which is a contradiction.

\vspace*{4mm}

\noindent So $v_{k+1}$ is a rich vertex verifying $N(v_{k+1}) \cap L_i(F) = N^-(v_{k+1}) \cap L_i(F)$ for all $i \leq k - 1$. Let $F_1 = F + yv_{k+1} + uv - vv_{k+1} - xv$ where $x$ is the predecessor of $v$ in $F$ and $y \in N^-(v_{k+1}) \cap U_{k-1}$ and let $F'$ be a maximal closure of $F_1$. Since $uk(F)$ is minimal, then $l_{F'} (z) = l_F (z)$ if $l_F (z) \leq k - 1$. This proves that $L_k(F') = (L_k(F)\setminus\{v_k\}) \cup {vk+1}$ and $v$ is still rich in $F'$ with $l_{F'} (v) \geq p \geq n$, but $v$ is an in-neighbor of $v_{k+1}$, then $\overline{v_{k+1}} = v$ and $h(v_{k+1}) \geq h(v)$. By supposing that $F$ is chosen such that $\sum_{w\;is\;bad} h(w)$ is maximal, we get $h(v_{k+1}) = h(v)$. This gives $l_{F'} (v) = l_{F_1} (v) = p$. Another important fact can be easily verified is that $l_{F'} (vs+1) = s$ for $k \leq s \leq p-1$. Hence $C_{v{k+1}} = C_v$. We repeat the same reasoning as above to prove that $v_{k'}$ ($k\leq k'\leq p$) is also a rich vertex verifying $N(v_{k'}) \cap L_i(F) = N^-(v_{k'}) \cap U_i$ for all $i \leq k - 1$. This can be verified by a simple induction for all the vertices in $C_v$.

\vspace*{6mm}

\noindent Step 3: If $\chi(D[C_v]) \geq l+2$, then by \emph{theorem 3.7} $D[C_v]$ contains a path $P(l, 1)$. This path can be completed to obtain a path $P(k, l)$ by adding $T_{v''}\cup v''v'$, where $v'$ is the end-vertex of the $P(l,1)$ corresponding to the block of length 1 which is rich and $v''$ is an in-neighbor of $v'$ of level $k-1$. Then we conclude that $\chi(D(Cv)) \leq l+1$.Color $C_v$ by the $l + 1$ colors $\{k,..., k + l\}$.

\vspace*{4mm}

\noindent If $C_v$ contains exactly $l + 2$ vertices (i.e. $f=1$), then at least two of the vertices of $C_v$ are not adjacent, we may suppose without loss of generality that $vv_j \notin E(G(D))$ since any vertex of $C_v$ can take the level $k$ in some convenient maximal spanning out-forest of $D$. We give each vertex $v_s$, $s \neq k$, the color $s$. Let $x \neq v_j$ be a rich vertex in $N(v) \cap U_j$ then we must have $x \in N^+(v) \cap L_j(F)$, otherwise we would use $x$ (as above) to make a directed path of length $k$ ending at $v$, and intersecting $C_v$ only at $v$, so by adding an appropriate directed path of $C_v$ we get a $P(k,l)$. We prove as above (as in (3)) that if $zw \in E(G(D))$ with $w \in T_v - T_x$ and $z \in T_x$, we have $z = x$, $w\rightarrow z$ and $l_F (w) < j$.Color a vertex $z \in T_x \cap U_i$ by the color $i + 1$ if $i < n - 1$ and by the color $j$ if $i = n-1$. We do the same for all rich vertices in $N(v) \cap U_j$ and the other non-reach vertices in $N(v) \cap U_j$ are colored by appropriate colors from $\{1,..., k - 1\}$ . The vertex $v$ is colored by $j$ and each remainder vertex $z \in U_i$ is colored by the color $i$, $k + 1 \leq i \leq k + l$. We get a good $(n - 1)$ coloring, which is a contradiction.

\vspace*{4mm}

\noindent We conclude that $l(C_v) > l+2$ (i.e. $f>1$), so $p > n$ and $l(C_V)=1+f(l+1)\geq 1+2(l+1)$. If we consider two vertices $v_s$ and $v_t$ in $C_v$ with $s < t \leq p$. Since $l(C_V)=\geq 1+2(l+1)$, then $C_v$ may be viewed as the union of two directed paths $Q_{v_sv_t}$ and $Q_{v_tv_s}$, such that one of them, say $P$, is of length at least $l + 1$. Set $S_{v_j} = T_{v_j} - T_{v_{j+1}}$ for $k \leq j \leq p - 1$ and $S_{v_p} = T_{v_p}$. If $x \in S_{v_t}$ and $y \in S_{v_s}$ are such that $xy \in E(G(D))$ and $\{x, y\} \neq \{v_s, v_t\}$. If $s\neq k$ or $y\neq v$, $P \cup P_{v_tx} \cup xy \cup P_w \cup wv_s \cup P_{v_sy}$ would contain a path $P(k, l)$ regardless of the orientation of $xy$, where $w$ is the in-neighbor of $v_s$ in $L_{k-1}$. So we must have $s = k$ and $y = v$. If $t \neq p$, $P \cup P_{v_tx} \cup xy\cup P_{z} \cup zu \cup uv \cup P_{vy} $ would contain a $P(k,l)$ regardless of the orientation of $xy$, where $z$ is the in-neighbor of $u=v_p$ in $L_{k-1}$. So we must have $t=p$ and $y=v$.

\vspace*{4mm}

\noindent Color $C_v$ by the $l+1$ colors $\{k,...,k+l\}$ such that $v$ is colored $k$ and $u$ is colored $k+1$. For all $w \in C_v$ of color $j = k + r$ we color each vertex $x \in L_m(S_w)$ by the color $k + h$ with $h \leq l$ and $h \equiv m+r-1$ mod $(l+1)$. We claim that the vertices in $S_u$ of color $k$ cannot be adjacent to $v$: If $w\in S_u$ is of color $k$ then $l_F(w)\geq p+l$ and so if $w$ is adjacent to $v$, then $P_v \cup vw \cup P_{vw}$ contain a $P(k,l)$ if $v\rightarrow w$ and $P_{vu} \cup uv \cup P_{uw} \cup wv$ contain a $P(k,l)$ if $w\rightarrow v$. Then this coloring is a good $n-1$-coloring of $D$, which contradicts the fact that $\chi(D) = n$. \hfill $\boxempty$

\section{Conclusion}

\noindent We have presented in this chapter the problem of finding paths with two blocks in $n$-chromatic digraphs. We have proven with two methods that for $n\geq4$, we can find any oriented path of length $n-1$ with two blocks in any $n$-chromatic digraph. What if we have more than two blocks?

\vspace*{4mm}

\noindent We conclude this chapter by stating this new conjecture of El-Sahili:

\vspace*{6mm}

\noindent {\bf Conjecture 3.8} ~\cite{pa3}: \emph{For $n \geq 8$, every $n$-chromatic digraph contains any oriented path of length $n - 1$.}

\vspace*{6mm}

\noindent In fact this conjecture generalizes Rosenfeld's conjecture which states that every tournament of order $n$ contains any oriented path of order $n-1$, which was proved by Havet and Thomass\'{e} with three exceptions which are tournaments of order 3,5 and 7. The condition $n\geq8$ is therefore necessary due to these three exceptions.

\end{document}